\documentclass[manuscript,printscheme,reqno]{aomart}
\usepackage{mystyle}
\usepackage{mathtools}
\usepackage{amssymb}
\usepackage{stmaryrd}
\usepackage{letltxmacro}

\newtheorem{thm}{Theorem}[section]

\newtheorem{lem}[thm]{Lemma}

\theoremstyle{definition}
\newtheorem{defn}[thm]{Definition}
\newtheorem{rem}[thm]{Remark}
\LetLtxMacro\origproof\proof
\LetLtxMacro\endorigproof\endproof
\makeatletter
\renewenvironment{proof}[1][\proofname]{\origproof[#1]\mbox{}\par\nobreak\@afterindentfalse\@afterheading\ignorespaces}{\endorigproof}
\makeatother

\DeclarePairedDelimiter{\floor}{\lfloor}{\rfloor}
\DeclarePairedDelimiter{\abs}{\lvert}{\rvert}
\DeclareMathOperator{\lcm}{lcm}
\newcommand{\N}{\mathbb{N}}
\newcommand{\Z}{\mathbb{Z}}
\newcommand{\R}{\mathbb{R}}

\allowdisplaybreaks
\title[Denumerant III]{Arithmetic Upper and Lower Bounds for the Denumerant Function III}
\author{Gerardo Ballesio}
\email{g.ballesio@gmail.com}
\givenname{Gerardo}
\surname{Ballesio}
\copyrightnote{}
\subject{primary}{matsc2020}{11P81}
\subject{secondary}{matsc2020}{05A17}
\subject{secondary}{matsc2020}{11N37}

\begin{document}

\begin{abstract}
	We extend the denumerant to a weighted denumerant, the sum of a function $\mu$ over the
	solutions in non-negative integers of $a_1x_1+\dotsb+a_kx_k=n$. For the monomial weights
	$x_1^{\alpha_1}\dotsm x_k^{\alpha_k}$, with exponents $\alpha_i\in\R_{\geq0}$, we prove upper and
	lower bounds of the kind of those of the first paper of this series (arXiv:2204.13689), with
	an explicit constant and explicit shifts. As an example of application, summing these bounds
	over all the numbers $a_1,\dotsc,a_k$ we find, by elementary means, the asymptotic main term of
	the sums $\sum_{y_1+\dotsb+y_k=n}\sigma_{\alpha_1}(y_1)\dotsm\sigma_{\alpha_k}(y_k)$ of products of
	divisor functions, for every $k$ and for real exponents $\alpha_i\neq0$ of either sign.
\end{abstract}
\maketitle
\setcounter{tocdepth}{2}
\tableofcontents

% =====================================================================
\section*{Notation}
% =====================================================================
Throughout the paper we use the following notation.
\begin{itemize}
	\item $\N=\{0,1,2,\dotsc\}$ and, for $t\in\N$, $\N_{\geq t}=\{t,t+1,t+2,\dotsc\}$; $\Z$ and $\R$
	      are the integers and the reals.
	\item $\floor{x}$ is the largest integer not exceeding $x\in\R$.
	\item $\delta_{x,y}$ is the Kronecker delta: $\delta_{x,y}=1$ if $x=y$ and $\delta_{x,y}=0$ if $x\neq y$.
	\item Empty sums are $0$ and empty products are $1$; in particular $0^0=1$.
	\item For $t\in\R$, $\R_{\geq t}$ and $\R_{>t}$ are the sets of the real numbers $y$ with $y\geq t$ and with $y>t$, respectively; in particular $\R_{\geq0}$ is the set of the non-negative real numbers.
	\item For $k\in\N_{\geq1}$, $x=(x_1,\dotsc,x_k)\in\N^k$ and
	      $\alpha=(\alpha_1,\dotsc,\alpha_k)\in\R^k$, we write $x^{\alpha}\coloneqq x_1^{\alpha_1}\dotsm x_k^{\alpha_k}$ and
	      $\abs{\alpha}\coloneqq\alpha_1+\dotsb+\alpha_k$; a negative exponent $\alpha_i$ is used only
	      with $x_i\geq1$.
	\item $\Gamma$ is Euler's Gamma function; recall that $\Gamma(\beta+1)=\beta!$ for $\beta\in\N$.
	\item For $\ell\in\N$ and $\beta_1,\dotsc,\beta_\ell\in\R_{\geq0}$, the multinomial coefficient is
	      \begin{equation*}
	      	\binom{\beta_1+\dotsb+\beta_\ell}{\beta_1,\dotsc,\beta_\ell}\coloneqq
	      	\frac{\Gamma(\beta_1+\dotsb+\beta_\ell+1)}{\Gamma(\beta_1+1)\dotsm\Gamma(\beta_\ell+1)},
	      \end{equation*}
	      which is $\frac{(\beta_1+\dotsb+\beta_\ell)!}{\beta_1!\dotsm\beta_\ell!}$ if
	      $\beta_1,\dotsc,\beta_\ell\in\N$.
	\item A letter without subscript, such as $x$, $\alpha$ or $a$, denotes a vector, and the same
	      letter with a subscript denotes one of its components: $x=(x_1,\dotsc,x_k)$; moreover
	      $0=(0,\dotsc,0)$.
	\item $a_1,a_2,a_3,\dotsc$ is a fixed sequence of positive integers.
\end{itemize}

The symbols that keep their meaning beyond a single proof are listed below, with the place where
they are defined. The other letters are local to the statement or the proof in which they appear.
\begin{center}
\renewcommand{\arraystretch}{1.25}
\begin{tabular}{@{}l p{0.52\textwidth} l@{}}
	\hline
	Symbol & Meaning & Defined in\\
	\hline
	$R_k(n)$ & solutions $x\in\N^k$ of $a_1x_1+\dotsb+a_kx_k=n$ & Section~\ref{s:intro}\\
	$D_k(n)$, $D_k^{\alpha}(n)$ & denumerant, and denumerant with the weight $x^{\alpha}$ & Section~\ref{s:intro}\\
	$d_k=d_k(a)$, $b_i$ & $\gcd(a_1,\dotsc,a_k)$ and $\lcm(d_{i-1},a_i)$ & Section~\ref{s:main}\\
	$\Delta_k$, $\nabla_k$ & the shifts of~\cite{ballesio-I} & Section~\ref{s:main}\\
	$\omega_k$, $E_k$ & the shift due to the weight, and the exponent & Section~\ref{s:main}\\
	$C_k^{\alpha}$ & the constant of the bounds & Section~\ref{s:main}\\
	$\alpha^{+}$, $\alpha^{-}$ & the positive and the negative part of $\alpha$ & Section~\ref{ss:setting}\\
	$E^{+}$, $\omega_k^{+}$, $C_{a}^{+}$, $D_{a}^{+}$ & the objects of Section~\ref{s:main} for the sequence $\alpha^{+}$ & Section~\ref{ss:setting}\\
	$B_{\beta,\gamma}$ & Euler's Beta function & Lemma~\ref{lem:beta}\\
	$\sigma_{u}(n)$, $\zeta(s)$ & sum of the $u$-th powers of the divisors, Riemann's zeta & Section~\ref{s:divisor}\\
	$S_{\alpha}(n)$ & the sum of products of divisor functions & Section~\ref{s:divisor}\\
	$\lVert\alpha\rVert$, $F$, $t_i$ & numbers that depend on $\alpha$ & Section~\ref{ss:setting}\\
	$H_{a}$, $h_{a}$ & the two shifts of Section~\ref{s:divisor} & Section~\ref{ss:setting}\\
	$c_{\alpha}$, $\kappa_{\alpha}$ & the constants of the main term & Lemma~\ref{lem:constant-sum}\\
	$s_i(\theta)$, $s_i'(\theta)$ & the exponents of the error terms & Theorem~\ref{thm:divisor}\\
	\hline
\end{tabular}
\end{center}

% =====================================================================
\section{The denumerant with a weight}\label{s:intro}
% =====================================================================
For $k\in\N_{\geq1}$ and $n\in\Z$ let
\begin{equation*}
	R_k(n)\coloneqq\bigl\{x\in\N^k : a_1x_1+\dotsb+a_kx_k=n\bigr\}
\end{equation*}
be the set of the solutions of $a_1x_1+\dotsb+a_kx_k=n$ in non-negative integers. It is finite,
since every solution satisfies $0\leq x_i\leq\frac{n}{a_i}$; and it is empty if $n<0$, since the
left side of the equation is non-negative. The denumerant is its number of
elements, $D_k(n)=\abs{R_k(n)}$, where each solution counts $1$. Here we give each solution a weight,
which may be any real number.

\begin{defn}\label{dfn:weighted}
	Let $k\in\N_{\geq1}$ and $\mu\colon\N^k\to\R$. For $n\in\Z$, the \emph{denumerant with weight
	$\mu$} is
	\begin{equation*}
		D_k^{\mu}(n)\coloneqq\sum_{x\in R_k(n)}\mu(x).
	\end{equation*}
\end{defn}

If $n<0$, the set $R_k(n)$ is empty, so the sum is empty and $D_k^{\mu}(n)=0$.

For $\mu=1$, every solution has weight $1$, and $D_k^{1}(n)=D_k(n)$. We consider weights given by
polynomials: a polynomial $P\in\R[X_1,\dotsc,X_k]$ defines the weight $x\mapsto P(x)$, and we
write $D_k^{P}(n)=\sum_{x\in R_k(n)}P(x)$ for the corresponding denumerant.

The simplest polynomials are the monomials $X^{\alpha}=X_1^{\alpha_1}\dotsm X_k^{\alpha_k}$, where
$X=(X_1,\dotsc,X_k)$ and $\alpha\in\N^k$. For them we write
\begin{equation*}
	D_k^{\alpha}(n)\coloneqq D_k^{X^{\alpha}}(n)=\sum_{x\in R_k(n)}x^{\alpha} .
\end{equation*}
For $\alpha=(0,\dotsc,0)=0$ we have $x^{\alpha}=1$, so $D_k^{0}(n)=D_k(n)$. Every polynomial is a
linear combination of monomials, and the next lemma reduces the denumerant with a polynomial
weight to the denumerants with monomial weights.

\begin{lem}\label{lem:linear}
	Let $k\in\N_{\geq1}$ and $P=\sum_{\alpha\in\Lambda}p_{\alpha} X^{\alpha}$, where $\Lambda\subseteq\N^k$ is finite
	and $p_{\alpha}\in\R$. Then, for every $n\in\Z$,
	\begin{equation*}
		D_k^{P}(n)=\sum_{\alpha\in\Lambda}p_{\alpha}\,D_k^{\alpha}(n).
	\end{equation*}
\end{lem}

\begin{proof}
	By the definition of $D_k^{P}$, then exchanging the two finite sums, and finally by the
	definition of $D_k^{\alpha}$,
	\begin{equation*}
		D_k^{P}(n)
		=\sum_{x\in R_k(n)}\sum_{\alpha\in\Lambda}p_{\alpha} x^{\alpha}
		=\sum_{\alpha\in\Lambda}p_{\alpha}\sum_{x\in R_k(n)}x^{\alpha}
		=\sum_{\alpha\in\Lambda}p_{\alpha}\,D_k^{\alpha}(n).\qedhere
	\end{equation*}
\end{proof}

So it is enough to study $D_k^{\alpha}(n)$ for $\alpha\in\N^k$. If the coefficients $p_{\alpha}$ are
non-negative, upper and lower bounds for the numbers $D_k^{\alpha}(n)$ give, by
Lemma~\ref{lem:linear}, upper and lower bounds for $D_k^{P}(n)$; if some of them are negative, an
upper bound for $D_k^{P}(n)$ takes the lower bounds of the corresponding $D_k^{\alpha}(n)$, and
conversely.

The sum that defines $D_k^{\alpha}(n)$ makes sense also when the exponents are not integers. So
from now on we let $\alpha\in\R_{\geq0}^k$, and
\begin{equation*}
	D_k^{\alpha}(n)\coloneqq\sum_{x\in R_k(n)}x_1^{\alpha_1}\dotsm x_k^{\alpha_k}
\end{equation*}
is the denumerant with the weight $x\mapsto x^{\alpha}$, a product of powers with
non-negative real exponents. The monomials are the case $\alpha\in\N^k$. The main result is
proved for real exponents, since the proof is the same.

As in the first paper of this series~\cite{ballesio-I}, we prove an upper and a lower bound with
a power (Section~\ref{s:main}, Theorem~\ref{thm:main}). As an example of application, in
Section~\ref{s:divisor} we sum these bounds over all the numbers $a_1,\dotsc,a_k$ and find the
asymptotic main term of the sums of products of divisor functions, also with negative exponents.
In the second paper of this series~\cite{ballesio-II} the bounds of the first paper are
sharpened with binomial coefficients in place of powers; here we use only those of the first
paper.

% =====================================================================
\section{Monomial weights}\label{s:main}
% =====================================================================
We use the integers of the first paper of this series~\cite{ballesio-I}. For $k\in\N_{\geq1}$ let
$d_k=d_k(a)\coloneqq\gcd(a_1,\dotsc,a_k)$, so that $d_1=a_1$, and $\Delta_k\coloneqq\sum_{i=2}^{k}(b_i-a_i)$,
where, for $i\in\N_{\geq2}$, $b_i\coloneqq\lcm(d_{i-1},a_i)$; in particular $\Delta_1=0$. The
notation $d_k(a)$, with $a=(a_1,\dotsc,a_k)$, recalls the dependence on the numbers $a_i$; we use it
only when this dependence has to be shown, as in Section~\ref{s:divisor}, where the $a_i$ vary, and
otherwise we write simply $d_k$. For
$k\in\N_{\geq2}$ let $\nabla_k\coloneqq b_2+\frac12\sum_{i=3}^{k}b_i$.

Like the numbers $a_i$, the exponents form a fixed sequence $\alpha=(\alpha_1,\alpha_2,\dotsc)$ of
numbers in $\R_{\geq0}$, and an object with the index $k$ uses only $\alpha_1,\dotsc,\alpha_k$; in
particular $\alpha$ may as well be a vector with at least $k$ components. For $k\in\N_{\geq1}$ we
set
\begin{equation*}
	E_k\coloneqq\alpha_1+\dotsb+\alpha_k+k-1,\qquad
	\omega_k\coloneqq\sum_{i=2}^{k}\bigl(1-\delta_{\alpha_i,0}\bigr)\,a_i,\qquad
	C_k^{\alpha}\coloneqq\frac{d_k}{(k-1)!\,
	a_1^{\alpha_1+1}\dotsm a_k^{\alpha_k+1}\dbinom{E_k}{\alpha_1,\dotsc,\alpha_k,k-1}}.
\end{equation*}
Like $\Delta_k$ and $\nabla_k$, the numbers $E_k$ and $\omega_k$ do not show the sequence on which
they depend; when we need to show it, as in Section~\ref{s:divisor}, we write $E_k^{\alpha}$ and
$\omega_k^{\alpha}$. Adding a variable, they change as follows:
\begin{equation}\label{eq:sigma-rec}
	E_{k+1}=E_k+\alpha_{k+1}+1,\qquad
	\omega_{k+1}=\omega_k+\bigl(1-\delta_{\alpha_{k+1},0}\bigr)\,a_{k+1},
\end{equation}
and $E_1=\alpha_1$, $\omega_1=0$.
For $\alpha=0$ the constant is $C_k^{0}=\frac{d_k}{(k-1)!\,a_1\dotsm a_k}$, the
constant of the denumerant.

\subsection{Preliminaries}\label{ss:mono-prel}
In this subsection we prepare the proof of Theorem~\ref{thm:main}, the polynomial bounds of
Subsection~\ref{ss:mono-poly}, which is by induction on $k$, as in the first paper.
The recursion for the denumerant, written with the weight of the last variable, reduces
$D_{k+1}^{\alpha}$ to a sum of values of $D_k^{\alpha}$ along an arithmetic
progression; after the inductive hypothesis, we are left with sums of the form
$\sum_u u^{\beta}(z-u)^{\gamma}$, which we estimate by comparison with an integral.

\subsubsection{The integral}
The integral that we compare with is Euler's Beta integral.

\begin{lem}\label{lem:beta}
	For $\beta,\gamma\in\R_{\geq0}$ let
	$B_{\beta,\gamma}\coloneqq\frac{\Gamma(\beta+1)\,\Gamma(\gamma+1)}{\Gamma(\beta+\gamma+2)}$, which is
	$\frac{\beta!\,\gamma!}{(\beta+\gamma+1)!}$ if $\beta,\gamma\in\N$. Then, for every $Y\in\R_{\geq0}$,
	\begin{equation*}
		\int_0^{Y}v^\beta(Y-v)^\gamma\,dv=B_{\beta,\gamma}\,Y^{\beta+\gamma+1}.
	\end{equation*}
\end{lem}

\begin{proof}
	For $Y=0$ both sides are $0$. Let $Y>0$. Putting $v=Y\nu$, so that $Y-v=Y(1-\nu)$ and $dv=Y\,d\nu$,
	\begin{equation*}
		\int_0^{Y}v^\beta(Y-v)^\gamma\,dv
		=Y^{\beta+\gamma+1}\int_0^{1}\nu^\beta(1-\nu)^\gamma\,d\nu,
	\end{equation*}
	and the last integral is Euler's Beta integral, which is equal to
	$\frac{\Gamma(\beta+1)\,\Gamma(\gamma+1)}{\Gamma(\beta+\gamma+2)}$ (see~\cite[Theorem~1.1.4]{andrews-askey-roy}).
	
	For integer exponents, which is the case of the polynomial weights, the lemma has an elementary
	proof, by parts, which we give for completeness. Let $\beta,\gamma\in\N$. We prove the statement, for
	every $\gamma\in\N$, by induction on $\beta$.
	
	\emph{Base case $\beta=0$.} The statement reads
	$\int_0^{Y}(Y-v)^\gamma\,dv=\frac{Y^{\gamma+1}}{\gamma+1}$, since $B_{0,\gamma}=\frac{\gamma!}{(\gamma+1)!}=\frac{1}{\gamma+1}$; and
	this holds, because $-\frac{(Y-v)^{\gamma+1}}{\gamma+1}$ is a primitive of $(Y-v)^\gamma$.
	
	\emph{Inductive step.} Let $\beta\in\N$, and assume that, for every $\gamma\in\N$,
	\begin{equation}\label{eq:beta-hyp}
		\int_0^{Y}v^\beta(Y-v)^\gamma\,dv=B_{\beta,\gamma}\,Y^{\beta+\gamma+1}.
	\end{equation}
	We have to prove that, for every $\gamma\in\N$,
	\begin{equation}\label{eq:beta-goal}
		\int_0^{Y}v^{\beta+1}(Y-v)^\gamma\,dv=B_{\beta+1,\gamma}\,Y^{\beta+\gamma+2}.
	\end{equation}
	We integrate by parts, differentiating $v^{\beta+1}$ and taking the primitive
	$-\frac{(Y-v)^{\gamma+1}}{\gamma+1}$ of $(Y-v)^\gamma$. The product $v^{\beta+1}(Y-v)^{\gamma+1}$ vanishes both at $v=0$
	and at $v=Y$, so the boundary term is $0$. Then we use Equation~\eqref{eq:beta-hyp} with $\gamma+1$ in
	place of $\gamma$, and finally the definition of $B_{\beta,\gamma+1}$:
	\begin{align*}
		\int_0^{Y}v^{\beta+1}(Y-v)^\gamma\,dv
		&=\frac{\beta+1}{\gamma+1}\int_0^{Y}v^{\beta}(Y-v)^{\gamma+1}\,dv
		=\frac{\beta+1}{\gamma+1}\,B_{\beta,\gamma+1}\,Y^{\beta+\gamma+2}\\
		&=\frac{\beta+1}{\gamma+1}\cdot\frac{\beta!\,(\gamma+1)!}{(\beta+\gamma+2)!}\,Y^{\beta+\gamma+2}\\
		&=\frac{(\beta+1)!\,\gamma!}{(\beta+\gamma+2)!}\,Y^{\beta+\gamma+2}
		=B_{\beta+1,\gamma}\,Y^{\beta+\gamma+2},
	\end{align*}
	which is Equation~\eqref{eq:beta-goal}.
\end{proof}

\subsubsection{An upper estimate for the sums}
In the upper bound the sum runs over points $u$ that differ by at least $1$. Each term is compared
with the integral over the interval of length $1$ centred at $u$, as in the midpoint rule; this
needs that one of the two factors $u^\beta$ and $(z-u)^\gamma$ is convex, and it is where the shift
$\frac12$ comes from. When neither is convex, we use the interval on the right of $u$, and the
shift is $1$.

\begin{lem}\label{lem:sum-upper}
	Let $\beta,\gamma\in\R_{\geq0}$, let $z\in\R_{\geq0}$, and let $U\subseteq[0,z]$ be a finite set of real
	numbers any two of which differ by at least $1$. Let $\eta\coloneqq\frac12$ if $\beta\in\R_{\geq1}$ or
	$\gamma\in\R_{\geq1}$, and let $\eta\coloneqq1$ otherwise. With $B_{\beta,\gamma}$ as in
	Lemma~\ref{lem:beta},
	\begin{equation*}
		\sum_{u\in U}u^\beta(z-u)^\gamma\leq B_{\beta,\gamma}\,(z+\eta)^{\beta+\gamma+1}.
	\end{equation*}
\end{lem}

\begin{proof}
	\begin{itemize}
		\item[$\scriptstyle\blacktriangleright$] \emph{Case~1: $\eta=1$.} The argument of this case does not use the size of $\beta$ and $\gamma$.
		Let $Y\coloneqq z+1$. By Lemma~\ref{lem:beta}, we have to
		prove that
		\begin{equation}\label{eq:sum-upper-right-goal}
			\sum_{u\in U}u^\beta(z-u)^\gamma\leq\int_0^{Y}v^\beta(Y-v)^\gamma\,dv .
		\end{equation}
		Let $u\in U$, and let $v$ be a point of $[u,u+1]$. Then:
		\begin{itemize}
			\item[$\bullet$] $v\geq u\geq0$, hence $v^\beta\geq u^\beta$;
			\item[$\bullet$] $v\leq u+1$, hence $Y-v\geq Y-u-1=z-u\geq0$, and $(Y-v)^\gamma\geq(z-u)^\gamma$.
		\end{itemize}
		Therefore $v^\beta(Y-v)^\gamma\geq u^\beta(z-u)^\gamma$ on $[u,u+1]$, an interval of length $1$, so
		\begin{equation}\label{eq:sum-upper-right-one}
			u^\beta(z-u)^\gamma\leq\int_{u}^{u+1}v^\beta(Y-v)^\gamma\,dv .
		\end{equation}
		Two points of $U$ differ by at least $1$, so the intervals $[u,u+1]$, for $u\in U$, overlap at
		most at an endpoint; and they lie in $[0,Y]$, since $0\leq u\leq z$. As the integrand is
		non-negative on $[0,Y]$, summing Equation~\eqref{eq:sum-upper-right-one} over $u\in U$ we get
		Equation~\eqref{eq:sum-upper-right-goal}.

		\item[$\scriptstyle\blacktriangleright$] \emph{Case~2: $\eta=\frac12$ and $\beta\in\R_{\geq1}$.} Let $Y\coloneqq z+\frac12$, and let $f$ be the function on
		$\bigl[-\frac12,Y\bigr]$ defined by
		\begin{equation*}
			f(v)\coloneqq v^\beta(Y-v)^\gamma\quad\text{for } 0\leq v\leq Y,
			\qquad
			f(v)\coloneqq0\quad\text{for } -\tfrac12\leq v<0 .
		\end{equation*}
		Then $f\geq0$, and by Lemma~\ref{lem:beta} its integral over $\bigl[-\frac12,Y\bigr]$ is
		$B_{\beta,\gamma}\,Y^{\beta+\gamma+1}$, the right-hand side of the statement. So we have to prove that
		\begin{equation}\label{eq:sum-upper-goal}
			\sum_{u\in U}u^\beta(z-u)^\gamma\leq\int_{-\frac12}^{Y}f(v)\,dv .
		\end{equation}
	
		\begin{itemize}
		\item[$\scriptstyle\blacklozenge$] \emph{Step~1: the values of $f$ at two symmetric points.} Let $u\in U$ and $0\leq \tau\leq\frac12$.
			We prove that
			\begin{equation}\label{eq:sum-upper-sym}
				f(u+\tau)+f(u-\tau)\geq2\,u^\beta(z-u)^\gamma .
			\end{equation}
			Since $0\leq u\leq z$ and $0\leq \tau\leq\frac12$, the point $u+\tau$ lies in $[0,Y]$, and
			\begin{equation*}
				Y-(u+\tau)=(z-u)+\Bigl(\frac12-\tau\Bigr)\geq z-u\geq0 .
			\end{equation*}
			Hence
			\begin{equation}\label{eq:sum-upper-plus}
				f(u+\tau)
				=(u+\tau)^\beta\bigl(Y-(u+\tau)\bigr)^\gamma
				\geq(u+\tau)^\beta(z-u)^\gamma .
			\end{equation}
			For $f(u-\tau)$ there are two cases.
			\begin{itemize}
				\item[$\bullet$] If $\tau\leq u$, then the point $u-\tau$ lies in $[0,Y]$, and
				\begin{equation*}
					Y-(u-\tau)=(z-u)+\frac12+\tau\geq z-u\geq0,
					\qquad\text{so}\qquad
					f(u-\tau)\geq(u-\tau)^\beta(z-u)^\gamma .
				\end{equation*}
				Since $\beta\in\R_{\geq1}$, the function $v\mapsto v^\beta$ is convex on $\R_{\geq0}$, so its
				value at $u$ is at most the average of its values at the two points $u+\tau$ and $u-\tau$,
				which are symmetric with respect to $u$:
				\begin{equation*}
					u^\beta\leq\frac{(u+\tau)^\beta+(u-\tau)^\beta}{2},
					\qquad\text{that is,}\qquad
					(u+\tau)^\beta+(u-\tau)^\beta\geq2\,u^\beta .
				\end{equation*}
				So, by Equation~\eqref{eq:sum-upper-plus},
				\begin{equation*}
					f(u+\tau)+f(u-\tau)
					\geq\bigl((u+\tau)^\beta+(u-\tau)^\beta\bigr)(z-u)^\gamma
					\geq2\,u^\beta(z-u)^\gamma .
				\end{equation*}
				\item[$\bullet$] If $\tau>u$, then $-\frac12\leq u-\tau<0$, since $u\geq0$ and $\tau\leq\frac12$; so $f(u-\tau)=0$. But now $u+\tau>2u$, and $2^\beta\geq2$ because
				$\beta\in\R_{\geq1}$. So, by Equation~\eqref{eq:sum-upper-plus},
				\begin{equation*}
					f(u+\tau)+f(u-\tau)
					\geq(u+\tau)^\beta(z-u)^\gamma
					\geq(2u)^\beta(z-u)^\gamma
					\geq2\,u^\beta(z-u)^\gamma .
				\end{equation*}
			\end{itemize}
			In both cases Equation~\eqref{eq:sum-upper-sym} holds.
	
		\item[$\scriptstyle\blacklozenge$] \emph{Step~2: one term and one interval.} Let $u\in U$. The interval $\bigl[u-\frac12,u+\frac12\bigr]$ lies in
			$\bigl[-\frac12,Y\bigr]$, where $f$ is defined, since $0\leq u\leq z$. We split it at $u$, we put $v=u+\tau$ in the right half and $v=u-\tau$ in the
			left half, and we use Equation~\eqref{eq:sum-upper-sym}:
			\begin{align}
				\int_{u-\frac12}^{u+\frac12}f(v)\,dv
				&=\int_{u}^{u+\frac12}f(v)\,dv+\int_{u-\frac12}^{u}f(v)\,dv
				=\int_0^{\frac12}f(u+\tau)\,d\tau+\int_0^{\frac12}f(u-\tau)\,d\tau\notag\\
				&=\int_0^{\frac12}\bigl(f(u+\tau)+f(u-\tau)\bigr)\,d\tau
				\geq\int_0^{\frac12}2\,u^\beta(z-u)^\gamma\,d\tau
				=u^\beta(z-u)^\gamma .\label{eq:sum-upper-one}
			\end{align}
			In the second integral, as $v$ goes from $u-\frac12$ to $u$, the variable $\tau=u-v$ goes from
			$\frac12$ to $0$, and $dv=-d\tau$; exchanging the two limits of integration removes the sign.
	
		\item[$\scriptstyle\blacklozenge$] \emph{Step~3: all the terms.} Two points of $U$ differ by at least $1$, so the intervals
			$\bigl[u-\frac12,u+\frac12\bigr]$, for $u\in U$, overlap at most at an endpoint; and they lie in
			$\bigl[-\frac12,Y\bigr]$, since $0\leq u\leq z$. As $f\geq0$, the sum of the integrals of $f$
			over these intervals is at most the integral of $f$ over $\bigl[-\frac12,Y\bigr]$. Hence, by
			Equation~\eqref{eq:sum-upper-one},
			\begin{equation*}
				\sum_{u\in U}u^\beta(z-u)^\gamma
				\leq\sum_{u\in U}\int_{u-\frac12}^{u+\frac12}f(v)\,dv
				\leq\int_{-\frac12}^{Y}f(v)\,dv,
			\end{equation*}
			which is Equation~\eqref{eq:sum-upper-goal}.
		\end{itemize}
		\item[$\scriptstyle\blacktriangleright$] \emph{Case~3: $\eta=\frac12$ and $\gamma\in\R_{\geq1}$.} We have to prove that
		\begin{equation}\label{eq:sum-upper-zero}
			\sum_{u\in U}u^\beta(z-u)^\gamma\leq B_{\beta,\gamma}\Bigl(z+\frac12\Bigr)^{\beta+\gamma+1}.
		\end{equation}
		The idea is to move the role of the weight from the factor $u^\beta$ to the other factor
		$(z-u)^\gamma$, which is convex since $\gamma\in\R_{\geq1}$. Let $U'$ be the set of the numbers $u'=z-u$ with
		$u\in U$. Then $U'\subseteq[0,z]$, and any two elements of $U'$ differ by at least $1$. For
		$u'=z-u$ we have $u=z-u'$, so $u^\beta(z-u)^\gamma=u'^\gamma(z-u')^\beta$. So we may apply Case~2 to
		the set $U'$, with the exponents $\gamma\in\R_{\geq1}$ and $\beta$ in place of $\beta$ and $\gamma$: in Case~2 the first
		exponent has to be at least $1$, and the second one is any non-negative real number. Moreover
		$B_{\gamma,\beta}=B_{\beta,\gamma}$, since the definition of $B_{\beta,\gamma}$ is symmetric in $\beta$ and $\gamma$. Hence
		\begin{equation*}
			\sum_{u\in U}u^\beta(z-u)^\gamma
			=\sum_{u'\in U'}u'^\gamma(z-u')^\beta
			\leq B_{\gamma,\beta}\Bigl(z+\frac12\Bigr)^{\beta+\gamma+1}
			=B_{\beta,\gamma}\Bigl(z+\frac12\Bigr)^{\beta+\gamma+1},
		\end{equation*}
		which is Equation~\eqref{eq:sum-upper-zero}.

	\end{itemize}
\end{proof}

\subsubsection{A lower estimate for the sums}
In the lower bound the points are $q+\varepsilon$, with $q\in\N$ and a fixed $\varepsilon\in[0,1)$. Each term is
compared with the integral over an interval of length at most $1$ on which the integrand is not
larger than the term.

\begin{lem}\label{lem:sum-lower}
	Let $\beta,\gamma\in\R_{\geq0}$, let $\varepsilon\in[0,1)$ and let $z\in\R_{\geq \varepsilon}$. With
	$B_{\beta,\gamma}$ as in Lemma~\ref{lem:beta}:
	\begin{enumerate}
		\item[\textup{(i)}] if $\beta=0$, then
		$\displaystyle\sum_{q=0}^{\floor{z-\varepsilon}}(z-\varepsilon-q)^\gamma\geq B_{0,\gamma}\,(z-\varepsilon)^{\gamma+1}$;
		\item[\textup{(ii)}] if $\beta\in\R_{>0}$ and $z\in\R_{\geq1}$, then
		$\displaystyle\sum_{q=0}^{\floor{z-\varepsilon}}(q+\varepsilon)^\beta(z-\varepsilon-q)^\gamma\geq B_{\beta,\gamma}\,(z-1)^{\beta+\gamma+1}$.
	\end{enumerate}
\end{lem}

\begin{proof}
	\emph{Part~\textup{(i)}.} Let $Y\coloneqq z-\varepsilon\geq0$. By Lemma~\ref{lem:beta} with $\beta=0$, we have to
	prove that
	\begin{equation}\label{eq:sum-lower-i-goal}
		\sum_{q=0}^{\floor{Y}}(Y-q)^\gamma\geq\int_0^{Y}(Y-v)^\gamma\,dv .
	\end{equation}
	We cut $[0,Y]$ at the integers: the intervals $[q,q+1]$ for $0\leq q<\floor{Y}$, and the last
	interval $[\floor{Y},Y]$. Each of them has length at most $1$, since $Y<\floor{Y}+1$. Let
	$0\leq q\leq\floor{Y}$, and let $v$ be a point of the interval that starts at $q$. Then
	$q\leq v\leq Y$, hence $0\leq Y-v\leq Y-q$ and $(Y-v)^\gamma\leq(Y-q)^\gamma$. So the integral of $(Y-v)^\gamma$
	over that interval is at most $(Y-q)^\gamma$ times its length, hence at most $(Y-q)^\gamma$. Summing over
	$0\leq q\leq\floor{Y}$, the integrals add up to the integral over $[0,Y]$, and we get
	Equation~\eqref{eq:sum-lower-i-goal}.
	
	\emph{Part~\textup{(ii)}.} Let $\beta\in\R_{>0}$ and $z\in\R_{\geq1}$, and let $u_q\coloneqq q+\varepsilon$, so that
	the sum of the statement is $\sum_{q=0}^{\floor{z-\varepsilon}}u_q^\beta(z-u_q)^\gamma$. By Lemma~\ref{lem:beta}
	with $Y=z-1\geq0$, we have to prove that
	\begin{equation}\label{eq:sum-lower-ii-goal}
		\sum_{q=0}^{\floor{z-\varepsilon}}u_q^\beta(z-u_q)^\gamma\geq\int_0^{z-1}v^\beta(z-1-v)^\gamma\,dv .
	\end{equation}
	
	\begin{itemize}
		\item[$\scriptstyle\blacklozenge$] \emph{Step~1: one term and one interval.} Let $q\in\N$ with
		$q\leq\floor{z-\varepsilon}$, so that $u_q\leq z$, and let $v$ be a point of $[u_q-1,u_q]$ that lies in $[0,z-1]$. Then:
		\begin{itemize}
			\item $0\leq v\leq u_q$, hence $v^\beta\leq u_q^\beta$;
			\item $v\geq u_q-1$, hence $z-1-v\leq z-u_q$; and $v\leq z-1$, hence $z-1-v\geq0$; so
			$(z-1-v)^\gamma\leq(z-u_q)^\gamma$.
		\end{itemize}
		Therefore $v^\beta(z-1-v)^\gamma\leq u_q^\beta(z-u_q)^\gamma$ at every such point $v$. These points form an
		interval of length at most $1$ (or no interval at all), so the integral of $v^\beta(z-1-v)^\gamma$ over
		the part of $[u_q-1,u_q]$ that lies in $[0,z-1]$ is at most $u_q^\beta(z-u_q)^\gamma$.
	
		\item[$\scriptstyle\blacklozenge$] \emph{Step~2: the intervals cover $[0,z-1]$.} The intervals $[u_q-1,u_q]$, for
		$0\leq q\leq\floor{z-\varepsilon}$, are consecutive, since $u_{q+1}-1=u_q$. The first one starts at
		$u_0-1=\varepsilon-1$, which is negative. The last one ends at $\floor{z-\varepsilon}+\varepsilon$, and
		\begin{equation*}
			\floor{z-\varepsilon}+\varepsilon>(z-\varepsilon-1)+\varepsilon=z-1 .
		\end{equation*}
		So together they cover $[0,z-1]$.
	
		\item[$\scriptstyle\blacklozenge$] \emph{Step~3: all the terms.} By Step~2, the integral over $[0,z-1]$ is the sum of the integrals
		over the parts of the intervals $[u_q-1,u_q]$ that lie in $[0,z-1]$; and by Step~1 each of these
		is at most $u_q^\beta(z-u_q)^\gamma$. This is Equation~\eqref{eq:sum-lower-ii-goal}.
	\end{itemize}
\end{proof}

In part~(ii) the shift is $1$ and not $\varepsilon$: the weight $u^\beta$ increases with $u$, so each term
controls the integral on the unit interval on its left, and not on its right.

\subsubsection{The recursion with a weight}
The next lemma is the
recursion of the first paper~\cite[Lemma~2.2]{ballesio-I}, where each term now carries the weight
of the last variable.

\begin{lem}\label{lem:recursion}
	Let $k\in\N_{\geq1}$, and let $n\in\N$ be a
	multiple of $d_{k+1}$. Let $r$ be the unique integer with $0\leq r<b_{k+1}$, $a_{k+1}\mid r$ and
	$d_k\mid n-r$, which exists by~\cite[Lemma~2.2]{ballesio-I} and satisfies
	$r\leq b_{k+1}-a_{k+1}$. Then
	\begin{equation*}
		D_{k+1}^{\alpha}(n)
		=\sum_{q=0}^{\floor*{\frac{n-r}{b_{k+1}}}}\biggl(\frac{r+qb_{k+1}}{a_{k+1}}\biggr)^{\alpha_{k+1}}
		D_k^{\alpha}(n-r-qb_{k+1}),
	\end{equation*}
	where every argument $n-r-qb_{k+1}$ is a non-negative multiple of $d_k$.
\end{lem}

\begin{proof}
	\begin{itemize}
		\item[$\scriptstyle\blacklozenge$] \emph{Step~1: the sum over $h$.} We classify the elements of $R_{k+1}(n)$ according to the
		value $h=a_{k+1}x_{k+1}$, which is a multiple of $a_{k+1}$ with $0\leq h\leq n$. For a fixed $h$, the last variable is
		$x_{k+1}=\frac{h}{a_{k+1}}$, and $(x_1,\dotsc,x_k)$ runs over $R_k(n-h)$; the weight of the
		solution is the product of $x_1^{\alpha_1}\dotsm x_k^{\alpha_k}$ and
		$\bigl(\frac{h}{a_{k+1}}\bigr)^{\alpha_{k+1}}$. Summing the weights
		first over $(x_1,\dotsc,x_k)\in R_k(n-h)$ and then over $h$,
		\begin{equation}\label{eq:recursion-h}
			D_{k+1}^{\alpha}(n)
			=\sum_{\substack{0\leq h\leq n\\ a_{k+1}\mid h}}\Bigl(\frac{h}{a_{k+1}}\Bigr)^{\alpha_{k+1}}D_k^{\alpha}(n-h).
		\end{equation}
		If $d_k\nmid n-h$, then $R_k(n-h)$ is empty and $D_k^{\alpha}(n-h)=0$; so only the $h$ with
		$d_k\mid n-h$ contribute.
	
		\item[$\scriptstyle\blacklozenge$] \emph{Step~2: the values of $h$ that contribute.} Let $h$ be a multiple of $a_{k+1}$. We have
		$d_k\mid n-r$ and $n-h=(n-r)-(h-r)$; hence $d_k\mid n-h$ if and only if $d_k\mid h-r$.
		Moreover
		$a_{k+1}\mid h-r$, because $a_{k+1}$ divides both $h$ and $r$. An integer is divisible by both
		$d_k$ and $a_{k+1}$ if and only if it is divisible by their least common multiple $b_{k+1}$.
		Hence $d_k\mid n-h$ if and only if $h=r+qb_{k+1}$ for some $q\in\Z$. Conversely, every
		number $r+qb_{k+1}$ is a multiple of $a_{k+1}$, since $a_{k+1}$ divides $r$ and $b_{k+1}$. Since
		$0\leq r<b_{k+1}$, the condition $0\leq r+qb_{k+1}\leq n$ reads $0\leq q\leq\floor*{\frac{n-r}{b_{k+1}}}$.
	
		\item[$\scriptstyle\blacklozenge$] \emph{Step~3: the formula.} By Steps~1 and~2, the terms of Equation~\eqref{eq:recursion-h} that
		contribute are those with $h=r+qb_{k+1}$ and $0\leq q\leq\floor*{\frac{n-r}{b_{k+1}}}$; this is
		the formula of the statement. For these $q$, the argument $n-h=n-r-qb_{k+1}$ is non-negative
		and a multiple of $d_k$.
	\end{itemize}
\end{proof}

\subsubsection{One step of the induction}
The upper and the lower bound of Theorem~\ref{thm:main} have the same form: the constant
$C_k^{\alpha}$ times a power of $n$ plus or minus a shift. In this section we show how a bound
of this form for $k$ gives one for $k+1$; the two lemmas are then applied in the proof of the
theorem, where only the shifts have to be followed. We first see how the constant changes.

\begin{lem}\label{lem:C-rec}
	Let $k\in\N_{\geq1}$. With $B_{\alpha_{k+1},E_k}$ as in Lemma~\ref{lem:beta},
	\begin{equation*}
		C_k^{\alpha}=\frac{d_k\,\Gamma(\alpha_1+1)\dotsm\Gamma(\alpha_k+1)}{\Gamma(E_k+1)\,a_1^{\alpha_1+1}\dotsm a_k^{\alpha_k+1}}
		\qquad\text{and}\qquad
		C_{k+1}^{\alpha}=C_k^{\alpha}\,\frac{B_{\alpha_{k+1},E_k}}{a_{k+1}^{\alpha_{k+1}}\,b_{k+1}} .
	\end{equation*}
\end{lem}

\begin{proof}
	\begin{itemize}
		\item[$\scriptstyle\blacklozenge$] \emph{Step~1: the multinomial coefficient.} We have $\alpha_1+\dotsb+\alpha_k+k=E_k+1$ and $\Gamma(k)=(k-1)!$.
		So, by the definition of the multinomial coefficient,
		\begin{align*}
			\binom{E_k}{\alpha_1,\dotsc,\alpha_k,k-1}(k-1)!
			&=\frac{\Gamma(E_k+1)}{\Gamma(\alpha_1+1)\dotsm\Gamma(\alpha_k+1)\,\Gamma(k)}\,(k-1)!\\
			&=\frac{\Gamma(E_k+1)}{\Gamma(\alpha_1+1)\dotsm\Gamma(\alpha_k+1)},
		\end{align*}
		and inserting this in the definition of $C_k^{\alpha}$ we get the first formula.
	
		\item[$\scriptstyle\blacklozenge$] \emph{Step~2: the recursion.} We have $\gcd(d_k,a_{k+1})=d_{k+1}$, hence
		$b_{k+1}=\lcm(d_k,a_{k+1})=\frac{d_k\,a_{k+1}}{d_{k+1}}$, that is,
		$\frac{d_{k+1}}{d_k\,a_{k+1}}=\frac{1}{b_{k+1}}$. Moreover $E_{k+1}=E_k+\alpha_{k+1}+1$ by Equation~\eqref{eq:sigma-rec}. So, by the
		first formula with $k+1$ in place of $k$, then
		separating the factors that contain $\alpha_{k+1}$, $a_{k+1}$ and $d_{k+1}$, and finally by the first
		formula, the definition of $B_{\alpha_{k+1},E_k}$ and the equation for $b_{k+1}$,
		\begin{align*}
			C_{k+1}^{\alpha}
			&=\frac{d_{k+1}\,\Gamma(\alpha_1+1)\dotsm\Gamma(\alpha_k+1)\,\Gamma(\alpha_{k+1}+1)}{\Gamma(E_{k+1}+1)\,a_1^{\alpha_1+1}\dotsm a_k^{\alpha_k+1}\,a_{k+1}^{\alpha_{k+1}+1}}\\
			&=\frac{d_k\,\Gamma(\alpha_1+1)\dotsm\Gamma(\alpha_k+1)}{\Gamma(E_k+1)\,a_1^{\alpha_1+1}\dotsm a_k^{\alpha_k+1}}
			\cdot\frac{\Gamma(\alpha_{k+1}+1)\,\Gamma(E_k+1)}{\Gamma(E_{k+1}+1)}\cdot\frac{d_{k+1}}{d_k\,a_{k+1}}\cdot\frac{1}{a_{k+1}^{\alpha_{k+1}}}\\
			&
			=C_k^{\alpha}\,\frac{B_{\alpha_{k+1},E_k}}{a_{k+1}^{\alpha_{k+1}}\,b_{k+1}} .\qedhere
		\end{align*}
	\end{itemize}
\end{proof}

\begin{lem}[one step, upper bound]\label{lem:step-upper}
	Let $k\in\N_{\geq1}$ and $N\in\R_{\geq0}$. Let $\eta\coloneqq\frac12$ if
	$\alpha_{k+1}\in\R_{\geq1}$ or $E_k\in\R_{\geq1}$, and let $\eta\coloneqq1$ otherwise. Assume that, for
	every multiple $m\in\N$ of $d_k$,
	\begin{equation}\label{eq:step-upper-hyp}
		D_k^{\alpha}(m)\leq C_k^{\alpha}\,(m+N)^{E_k}.
	\end{equation}
	Then, for every multiple $n\in\N$ of $d_{k+1}$,
	\begin{equation}\label{eq:step-upper-goal}
		D_{k+1}^{\alpha}(n)\leq C_{k+1}^{\alpha}(n+N+\eta\,b_{k+1})^{E_{k+1}}.
	\end{equation}
\end{lem}

\begin{proof}
	Let $n\in\N$ be a multiple of $d_{k+1}$, and let $r$ be as in Lemma~\ref{lem:recursion}. We have
	to prove Equation~\eqref{eq:step-upper-goal}. We set
	\begin{equation*}
		z\coloneqq\frac{n+N}{b_{k+1}},
		\qquad
		u_q\coloneqq q+\frac{r}{b_{k+1}}\quad\text{for } q\in\N .
	\end{equation*}
	
	\begin{itemize}
		\item[$\scriptstyle\blacklozenge$] \emph{Step~1: one term of the recursion.} Let $0\leq q\leq\floor*{\frac{n-r}{b_{k+1}}}$. By the
		definition of $u_q$, the weight of the term with index $q$ in Lemma~\ref{lem:recursion} is
		\begin{equation}\label{eq:step-upper-weight}
			\biggl(\frac{r+qb_{k+1}}{a_{k+1}}\biggr)^{\alpha_{k+1}}
			=\biggl(\frac{b_{k+1}}{a_{k+1}}\,u_q\biggr)^{\alpha_{k+1}}
			=\Bigl(\frac{b_{k+1}}{a_{k+1}}\Bigr)^{\alpha_{k+1}}u_q^{\alpha_{k+1}}.
		\end{equation}
		The argument $n-r-qb_{k+1}$ is a non-negative multiple of $d_k$ (Lemma~\ref{lem:recursion}), so
		we may apply Equation~\eqref{eq:step-upper-hyp} with $m=n-r-qb_{k+1}$; then, by the definitions
		of $z$ and $u_q$,
		\begin{equation}\label{eq:step-upper-term}
			D_k^{\alpha}(n-r-qb_{k+1})
			\leq C_k^{\alpha}\,(n-r-qb_{k+1}+N)^{E_k}
			=C_k^{\alpha}\,b_{k+1}^{E_k}\,(z-u_q)^{E_k}.
		\end{equation}
	
		\item[$\scriptstyle\blacklozenge$] \emph{Step~2: the sum.} By Lemma~\ref{lem:recursion}, and then by
		Equations~\eqref{eq:step-upper-weight} and~\eqref{eq:step-upper-term}, where the inequality
		is kept because the weights are non-negative,
		\begin{align}
			D_{k+1}^{\alpha}(n)
			&=\sum_{q=0}^{\floor*{\frac{n-r}{b_{k+1}}}}\biggl(\frac{r+qb_{k+1}}{a_{k+1}}\biggr)^{\alpha_{k+1}}
			D_k^{\alpha}(n-r-qb_{k+1})\notag\\
			&\leq C_k^{\alpha}\Bigl(\frac{b_{k+1}}{a_{k+1}}\Bigr)^{\alpha_{k+1}}b_{k+1}^{E_k}
			\sum_{q=0}^{\floor*{\frac{n-r}{b_{k+1}}}}u_q^{\alpha_{k+1}}(z-u_q)^{E_k}.\label{eq:step-upper-sum}
		\end{align}
	
		\item[$\scriptstyle\blacklozenge$] \emph{Step~3: the estimate of the sum.} Let $U$ be the set of the numbers $u_q$ with
		$0\leq q\leq\floor*{\frac{n-r}{b_{k+1}}}$. Two of them differ by at least $1$. Moreover
		$0\leq r+qb_{k+1}\leq n$ for these $q$, so $0\leq u_q\leq\frac{n}{b_{k+1}}\leq z$, because
		$N\geq0$. Hence $U\subseteq[0,z]$, and we may apply Lemma~\ref{lem:sum-upper} with $\beta=\alpha_{k+1}$ and $\gamma=E_k$;
		the number $\eta$ of that lemma is then the number $\eta$ of the present statement. It
		gives
		\begin{equation}\label{eq:step-upper-est}
			\sum_{q=0}^{\floor*{\frac{n-r}{b_{k+1}}}}u_q^{\alpha_{k+1}}(z-u_q)^{E_k}
			\leq B_{\alpha_{k+1},E_k}(z+\eta)^{E_{k+1}}.
		\end{equation}
	
		\item[$\scriptstyle\blacklozenge$] \emph{Step~4: conclusion.} By Equations~\eqref{eq:step-upper-sum} and~\eqref{eq:step-upper-est},
		then collecting the powers of $b_{k+1}$, where we use
		\begin{equation*}
			b_{k+1}^{\alpha_{k+1}+E_k}=\frac{1}{b_{k+1}}\,b_{k+1}^{E_{k+1}},
		\end{equation*}
		then by $b_{k+1}\,(z+\eta)=n+N+\eta\,b_{k+1}$, and finally by Lemma~\ref{lem:C-rec},
		\begin{align*}
			D_{k+1}^{\alpha}(n)
			&\leq C_k^{\alpha}\Bigl(\frac{b_{k+1}}{a_{k+1}}\Bigr)^{\alpha_{k+1}}b_{k+1}^{E_k}\,
			B_{\alpha_{k+1},E_k}(z+\eta)^{E_{k+1}}\\
			&=C_k^{\alpha}\,\frac{B_{\alpha_{k+1},E_k}}{a_{k+1}^{\alpha_{k+1}}\,b_{k+1}}
			\bigl(b_{k+1}\,(z+\eta)\bigr)^{E_{k+1}}\\
			&=C_k^{\alpha}\,\frac{B_{\alpha_{k+1},E_k}}{a_{k+1}^{\alpha_{k+1}}\,b_{k+1}}
			(n+N+\eta\,b_{k+1})^{E_{k+1}}\\
			&=C_{k+1}^{\alpha}(n+N+\eta\,b_{k+1})^{E_{k+1}},
		\end{align*}
		which is Equation~\eqref{eq:step-upper-goal}.
	\end{itemize}
\end{proof}

\begin{lem}[one step, lower bound]\label{lem:step-lower}
	Let $k\in\N_{\geq1}$ and $N\in\R_{\geq0}$. Assume that, for every multiple
	$m\in\N$ of $d_k$ with $m\geq N$,
	\begin{equation}\label{eq:step-lower-hyp}
		D_k^{\alpha}(m)\geq C_k^{\alpha}\,(m-N)^{E_k}.
	\end{equation}
	Let $M\coloneqq N+b_{k+1}-a_{k+1}$ if $\alpha_{k+1}=0$, and $M\coloneqq N+b_{k+1}$ if $\alpha_{k+1}\in\R_{>0}$. Then,
	for every multiple $n\in\N$ of $d_{k+1}$ with $n\geq M$,
	\begin{equation}\label{eq:step-lower-goal}
		D_{k+1}^{\alpha}(n)\geq C_{k+1}^{\alpha}\,(n-M)^{E_{k+1}}.
	\end{equation}
\end{lem}

\begin{proof}
	Let $n\geq M$ be a multiple of $d_{k+1}$, and let $r$ be as in Lemma~\ref{lem:recursion}. We
	have to prove Equation~\eqref{eq:step-lower-goal}. We set
	\begin{equation*}
		z\coloneqq\frac{n-N}{b_{k+1}},
		\qquad
		\varepsilon\coloneqq\frac{r}{b_{k+1}} .
	\end{equation*}
	Since $0\leq r\leq b_{k+1}-a_{k+1}$ (Lemma~\ref{lem:recursion}), we have $\varepsilon\in[0,1)$. Moreover,
	in both cases $M\geq N+b_{k+1}-a_{k+1}$, so
	\begin{equation}\label{eq:step-lower-z}
		z-\varepsilon
		=\frac{n-N-r}{b_{k+1}}
		\geq\frac{n-N-(b_{k+1}-a_{k+1})}{b_{k+1}}
		\geq\frac{n-M}{b_{k+1}}
		\geq0 .
	\end{equation}
	
	\begin{itemize}
		\item[$\scriptstyle\blacklozenge$] \emph{Step~1: the terms that we keep.} Let $0\leq q\leq\floor{z-\varepsilon}$. Then $q\leq z-\varepsilon$, and, by
		the definitions of $z$ and $\varepsilon$,
		\begin{equation}\label{eq:step-lower-arg}
			n-r-qb_{k+1}
			=N+b_{k+1}\,(z-\varepsilon-q)
			\geq N
			\geq0 .
		\end{equation}
		In particular $r+qb_{k+1}\leq n$, so $q\leq\floor*{\frac{n-r}{b_{k+1}}}$: the term with index
		$q$ belongs to the sum of Lemma~\ref{lem:recursion}. The other terms of that sum are
		non-negative, so we may leave them out:
		\begin{equation}\label{eq:step-lower-keep}
			D_{k+1}^{\alpha}(n)
			\geq\sum_{q=0}^{\floor{z-\varepsilon}}\biggl(\frac{r+qb_{k+1}}{a_{k+1}}\biggr)^{\alpha_{k+1}}
			D_k^{\alpha}(n-r-qb_{k+1}).
		\end{equation}
	
		\item[$\scriptstyle\blacklozenge$] \emph{Step~2: one term.} Let $0\leq q\leq\floor{z-\varepsilon}$. By the definition of $\varepsilon$, the weight is
		\begin{equation}\label{eq:step-lower-weight}
			\biggl(\frac{r+qb_{k+1}}{a_{k+1}}\biggr)^{\alpha_{k+1}}
			=\Bigl(\frac{b_{k+1}}{a_{k+1}}\Bigr)^{\alpha_{k+1}}(q+\varepsilon)^{\alpha_{k+1}}.
		\end{equation}
		By Equation~\eqref{eq:step-lower-arg}, the argument $m=n-r-qb_{k+1}$ is a multiple of $d_k$
		(Lemma~\ref{lem:recursion}) with $m\geq N$ and $m-N=b_{k+1}(z-\varepsilon-q)$; so
		Equation~\eqref{eq:step-lower-hyp} gives
		\begin{equation}\label{eq:step-lower-term}
			D_k^{\alpha}(n-r-qb_{k+1})
			\geq C_k^{\alpha}\,(n-r-qb_{k+1}-N)^{E_k}
			=C_k^{\alpha}\,b_{k+1}^{E_k}\,(z-\varepsilon-q)^{E_k}.
		\end{equation}
		Inserting Equations~\eqref{eq:step-lower-weight} and~\eqref{eq:step-lower-term} in
		Equation~\eqref{eq:step-lower-keep},
		\begin{equation}\label{eq:step-lower-sum}
			D_{k+1}^{\alpha}(n)
			\geq C_k^{\alpha}\Bigl(\frac{b_{k+1}}{a_{k+1}}\Bigr)^{\alpha_{k+1}}b_{k+1}^{E_k}
			\sum_{q=0}^{\floor{z-\varepsilon}}(q+\varepsilon)^{\alpha_{k+1}}(z-\varepsilon-q)^{E_k}.
		\end{equation}
	
		\item[$\scriptstyle\blacklozenge$] \emph{Step~3: the estimate of the sum.} We prove that, in both cases,
		\begin{equation}\label{eq:step-lower-est}
			\sum_{q=0}^{\floor{z-\varepsilon}}(q+\varepsilon)^{\alpha_{k+1}}(z-\varepsilon-q)^{E_k}
			\geq B_{\alpha_{k+1},E_k}\Bigl(\frac{n-M}{b_{k+1}}\Bigr)^{E_{k+1}}.
		\end{equation}
		\begin{itemize}
			\item If $\alpha_{k+1}=0$, then by Lemma~\ref{lem:sum-lower}~(i) with $\gamma=E_k$ the sum is at least
			$B_{0,E_k}\,(z-\varepsilon)^{E_k+1}$, where $E_k+1=E_{k+1}$ and $B_{0,E_k}=B_{\alpha_{k+1},E_k}$ since $\alpha_{k+1}=0$; and $z-\varepsilon\geq\frac{n-M}{b_{k+1}}\geq0$ by
			Equation~\eqref{eq:step-lower-z}.
			\item If $\alpha_{k+1}\in\R_{>0}$, then $M=N+b_{k+1}$, so $z-1=\frac{n-N-b_{k+1}}{b_{k+1}}=\frac{n-M}{b_{k+1}}\geq0$;
			in particular $z\geq1$. By Lemma~\ref{lem:sum-lower}~(ii) with $\beta=\alpha_{k+1}$ and $\gamma=E_k$, the sum is
			at least $B_{\alpha_{k+1},E_k}\,(z-1)^{E_{k+1}}$.
		\end{itemize}
	
		\item[$\scriptstyle\blacklozenge$] \emph{Step~4: conclusion.} By Equations~\eqref{eq:step-lower-sum} and~\eqref{eq:step-lower-est},
		then collecting the powers of $b_{k+1}$, and finally by Lemma~\ref{lem:C-rec},
		\begin{align*}
			D_{k+1}^{\alpha}(n)
			&\geq C_k^{\alpha}\Bigl(\frac{b_{k+1}}{a_{k+1}}\Bigr)^{\alpha_{k+1}}b_{k+1}^{E_k}\,
			B_{\alpha_{k+1},E_k}\Bigl(\frac{n-M}{b_{k+1}}\Bigr)^{E_{k+1}}\\
			&=C_k^{\alpha}\,\frac{B_{\alpha_{k+1},E_k}}{a_{k+1}^{\alpha_{k+1}}\,b_{k+1}}\,(n-M)^{E_{k+1}}
			=C_{k+1}^{\alpha}\,(n-M)^{E_{k+1}},
		\end{align*}
		which is Equation~\eqref{eq:step-lower-goal}.
	\end{itemize}
\end{proof}

\subsection{Polynomial bounds}\label{ss:mono-poly}
\begin{thm}\label{thm:main}
	Let $k\in\N_{\geq2}$, and let $n\in\N$ be a multiple of $d_k$. Then
	\begin{enumerate}
		\item[\textup{(i)}] $D_k^{\alpha}(n)\leq C_k^{\alpha}\,(n+\nabla_k)^{E_k}$;
		\item[\textup{(ii)}] if $n\geq\Delta_k+\omega_k$, then
		$D_k^{\alpha}(n)\geq C_k^{\alpha}\,(n-\Delta_k-\omega_k)^{E_k}$.
	\end{enumerate}
\end{thm}

For $\alpha=0$ we have $\omega_k=0$ and
$C_k^{0}=\frac{d_k}{(k-1)!\,a_1\dotsm a_k}=\frac{1}{(k-1)!\,b_2\dotsm b_k}$, since
$b_2\dotsm b_k=\frac{a_1\dotsm a_k}{d_k}$ (in $b_i=\frac{d_{i-1}a_i}{d_i}$ the quotients
$\frac{d_{i-1}}{d_i}$ cancel in a chain, and $d_1=a_1$); so Theorem~\ref{thm:main} gives the upper and the
lower bound of~\cite[Theorem~1.3]{ballesio-I}, the latter with $\geq$ in place of $>$. The shift
$\omega_k$ in the lower bound cannot be dropped: a weight $x_i^{\alpha_i}$ with
$\alpha_i\in\R_{>0}$ vanishes on the solutions with $x_i=0$. For instance, for $a_1=a_2=1$ and
$\alpha=(1,1)$ we have $\omega_2=1$ and
$D_2^{\alpha}(n)=\sum_{x=0}^{n}x(n-x)=\frac{n^3-n}{6}$, which is less than
$C_2^{\alpha}\,n^3=\frac{n^3}{6}$ for every $n\geq1$, while the bound of the theorem is
$\frac{(n-1)^3}{6}$. The first variable needs no shift, since
$x_1$ is determined by the others.

\subsubsection{Proof of Theorem~\ref{thm:main}}
The sequence $\alpha$ is fixed, and both parts are proved by induction on $k$. Both inductions
start from one variable, where the denumerant with a weight is known exactly. Let $m\in\N$ be a
multiple of $d_1=a_1$. The equation $a_1x_1=m$ has the single solution $x_1=\frac{m}{a_1}$, whose
weight is $\bigl(\frac{m}{a_1}\bigr)^{\alpha_1}$; and
$C_1^{\alpha}=\frac{d_1}{a_1^{\alpha_1+1}}=\frac{1}{a_1^{\alpha_1}}$ by Lemma~\ref{lem:C-rec},
since $E_1=\alpha_1$ and $d_1=a_1$. Hence
\begin{equation}\label{eq:k1}
	D_1^{\alpha}(m)
	=\Bigl(\frac{m}{a_1}\Bigr)^{\alpha_1}
	=C_1^{\alpha}\,m^{E_1}.
\end{equation}

\begin{proof}[Proof of Theorem~\ref{thm:main}~(i)]
	\emph{Base case $k=2$.} Let $n\in\N$ be a multiple of $d_2$. Since $\nabla_2=b_2$, we have to
	prove that
	\begin{equation}\label{eq:upper-base-goal}
		D_2^{\alpha}(n)\leq C_2^{\alpha}\,(n+b_2)^{E_2}.
	\end{equation}
	We apply Lemma~\ref{lem:step-upper} with $k=1$ and $N=0$. Its assumption,
	Equation~\eqref{eq:step-upper-hyp}, reads $D_1^{\alpha}(m)\leq C_1^{\alpha}\,m^{E_1}$, which
	holds by Equation~\eqref{eq:k1}. So the lemma gives the first inequality below, where $\eta$ is
	$\frac12$ or $1$; the second one holds because $\eta\leq1$:
	\begin{equation*}
		D_2^{\alpha}(n)
		\leq C_2^{\alpha}\,(n+\eta\,b_2)^{E_2}
		\leq C_2^{\alpha}\,(n+b_2)^{E_2},
	\end{equation*}
	which is Equation~\eqref{eq:upper-base-goal}.

	\emph{Inductive step.} Let $k\geq2$, and assume that, for every multiple $m\in\N$ of $d_k$,
	\begin{equation}\label{eq:upper-hyp}
		D_k^{\alpha}(m)\leq C_k^{\alpha}\,(m+\nabla_k)^{E_k}.
	\end{equation}
	We have to prove that, for every multiple $n\in\N$ of $d_{k+1}$,
	\begin{equation}\label{eq:upper-goal}
		D_{k+1}^{\alpha}(n)\leq C_{k+1}^{\alpha}\,(n+\nabla_{k+1})^{E_{k+1}}.
	\end{equation}
	We apply Lemma~\ref{lem:step-upper} with $N=\nabla_k\geq0$; its assumption is
	Equation~\eqref{eq:upper-hyp}. Here $E_k\geq k-1\geq1$, so the number $\eta$ of the lemma is
	$\frac12$, whatever $\alpha_{k+1}$ is. By the definition of $\nabla_k$, we have
	$\nabla_k+\frac12\,b_{k+1}=\nabla_{k+1}$. So the lemma gives
	\begin{equation*}
		D_{k+1}^{\alpha}(n)
		\leq C_{k+1}^{\alpha}\Bigl(n+\nabla_k+\frac12\,b_{k+1}\Bigr)^{E_{k+1}}
		=C_{k+1}^{\alpha}\,(n+\nabla_{k+1})^{E_{k+1}},
	\end{equation*}
	which is Equation~\eqref{eq:upper-goal}.
\end{proof}

\begin{proof}[Proof of Theorem~\ref{thm:main}~(ii)]
	In both steps we use Equation~\eqref{eq:sigma-rec}: adding the variable $x_{k+1}$, the shift
	$\omega_k$ does not change if $\alpha_{k+1}=0$, and it grows by $a_{k+1}$ if
	$\alpha_{k+1}\in\R_{>0}$.

	\emph{Base case $k=2$.} Let $n\in\N$ be a multiple of $d_2$ with $n\geq\Delta_2+\omega_2$. We
	have to prove that
	\begin{equation}\label{eq:lower-base-goal}
		D_2^{\alpha}(n)\geq C_2^{\alpha}\,(n-\Delta_2-\omega_2)^{E_2}.
	\end{equation}
	We apply Lemma~\ref{lem:step-lower} with $k=1$ and $N=0$. Its assumption,
	Equation~\eqref{eq:step-lower-hyp}, reads $D_1^{\alpha}(m)\geq C_1^{\alpha}\,m^{E_1}$, which
	holds by Equation~\eqref{eq:k1}. The number $M$ of the lemma is $b_2-a_2$ if $\alpha_2=0$ and
	$b_2$ if $\alpha_2\in\R_{>0}$. Since $\Delta_2=b_2-a_2$ and $\omega_1=0$, by
	Equation~\eqref{eq:sigma-rec} in both cases
	\begin{equation*}
		M=\Delta_2+\omega_2 .
	\end{equation*}
	So the conclusion of the lemma, Equation~\eqref{eq:step-lower-goal}, is
	Equation~\eqref{eq:lower-base-goal}.

	\emph{Inductive step.} Let $k\geq2$, and assume that, for every multiple $m\in\N$ of $d_k$ with
	$m\geq\Delta_k+\omega_k$,
	\begin{equation}\label{eq:lower-hyp}
		D_k^{\alpha}(m)\geq C_k^{\alpha}\,(m-\Delta_k-\omega_k)^{E_k}.
	\end{equation}
	We have to prove that, for every multiple $n\in\N$ of $d_{k+1}$ with
	$n\geq\Delta_{k+1}+\omega_{k+1}$,
	\begin{equation}\label{eq:lower-goal}
		D_{k+1}^{\alpha}(n)\geq C_{k+1}^{\alpha}\,(n-\Delta_{k+1}-\omega_{k+1})^{E_{k+1}}.
	\end{equation}
	We apply Lemma~\ref{lem:step-lower} with $N=\Delta_k+\omega_k\geq0$; its assumption is
	Equation~\eqref{eq:lower-hyp}. We compute the number $M$ of the lemma, using
	$\Delta_{k+1}=\Delta_k+b_{k+1}-a_{k+1}$ and Equation~\eqref{eq:sigma-rec}:
	\begin{itemize}
		\item if $\alpha_{k+1}=0$, then
		$M=\Delta_k+\omega_k+b_{k+1}-a_{k+1}=\Delta_{k+1}+\omega_{k+1}$;
		\item if $\alpha_{k+1}\in\R_{>0}$, then
		$M=\Delta_k+\omega_k+b_{k+1}=\Delta_{k+1}+\omega_k+a_{k+1}=\Delta_{k+1}+\omega_{k+1}$.
	\end{itemize}
	So the conclusion of the lemma, Equation~\eqref{eq:step-lower-goal}, is
	Equation~\eqref{eq:lower-goal}.
\end{proof}

% =====================================================================
\section{An application: the main term of sums of products of divisor functions}\label{s:divisor}
% =====================================================================
As an example of application of Theorem~\ref{thm:main}, we find the dominant term, as
$n\to\infty$, of a classical arithmetic sum. In this section the numbers $a_1,\dotsc,a_k$ are no
longer fixed: we sum the denumerant with a monomial weight over all of them. For $u\in\R$ and
$n\in\N_{\geq1}$ let $\sigma_{u}(n)\coloneqq\sum_{d\mid n}d^{u}$ be the sum of the $u$-th
powers of the divisors of $n$, and for $s\in\R_{>1}$ let $\zeta(s)\coloneqq\sum_{n\in\N_{\geq1}}n^{-s}$.
Throughout the section $k\in\N_{\geq2}$ and $\alpha\in(\R\setminus\{0\})^k$, that is, every
exponent $\alpha_i$ is a real number different from $0$, positive or negative. We study, for
$n\in\N_{\geq1}$,
\begin{equation*}
	S_{\alpha}(n)\coloneqq\sum_{\substack{y\in\N_{\geq1}^k\\y_1+\dotsb+y_k=n}}\sigma_{\alpha_1}(y_1)\dotsm\sigma_{\alpha_k}(y_k).
\end{equation*}
For $k=2$ Ramanujan~\cite{ramanujan} proved an asymptotic formula for $S_{\alpha}(n)$ when
$\alpha_1$ and $\alpha_2$ are odd positive integers; Ingham~\cite{ingham} proved its main term for
all positive real $\alpha_1,\alpha_2$; and Lemke Oliver, Shrestha and Thorne~\cite{lemke-oliver}
gave a new proof, which also obtains lower-order terms for most ranges of the parameters. Here we obtain the main term for every
$k$, with explicit bounds, by elementary means: from Theorem~\ref{thm:main} alone, without modular
forms or complex analysis.

\subsection{The setting}\label{ss:setting}
We split $\alpha$ into its positive and its negative part,
\begin{equation*}
	\alpha^{+}\coloneqq\bigl(\max\{\alpha_1,0\},\dotsc,\max\{\alpha_k,0\}\bigr),\qquad
	\alpha^{-}\coloneqq\bigl(\min\{\alpha_1,0\},\dotsc,\min\{\alpha_k,0\}\bigr),
\end{equation*}
so that $\alpha=\alpha^{+}+\alpha^{-}$. Since every $\alpha_i$ is different from $0$, we have
$\alpha_i^{+}=0$ exactly when $\alpha_i<0$. The objects of Section~\ref{s:main} will be used for the
sequence $\alpha^{+}$, with non-negative exponents, and we mark them with $+$:
\begin{equation*}
	E^{+}\coloneqq E_k^{\alpha^{+}}=\abs{\alpha^{+}}+k-1,\qquad
	\omega_k^{+}\coloneqq\omega_k^{\alpha^{+}},
\end{equation*}
Moreover, the numbers $a_1,\dotsc,a_k$ now vary, and we sum over them. So we write them as
a vector $a\in\N_{\geq1}^k$, with components $a_1,\dotsc,a_k$, and we add $a$ as an index to the two objects that
depend on them and appear in these sums:
\begin{equation*}
	C_{a}^{+}\coloneqq C_k^{\alpha^{+}},\qquad
	D_{a}^{+}\coloneqq D_k^{\alpha^{+}},
\end{equation*}
both for the numbers $a_1,\dotsc,a_k$; recall that $D_{a}^{+}(m)=0$ for $m\in\Z_{<0}$. The numbers $d_k=d_k(a)$ and $b_i$ are those of
Section~\ref{s:main}, built from the components $a_i$ of $a$. We also
use
\begin{gather*}
	\lVert\alpha\rVert\coloneqq\abs{\alpha_1}+\dotsb+\abs{\alpha_k},\qquad
	F\coloneqq\lVert\alpha\rVert+k-1,\qquad
	t_i\coloneqq\abs{\alpha_i}+1\quad\text{for }1\leq i\leq k,\\
	H_{a}\coloneqq a_1+b_2+\dotsb+b_k,\qquad
	h_{a}\coloneqq\sum_{i\in\{1,\dotsc,k\}\colon\alpha_i<0}a_i .
\end{gather*}
Every $t_i$ is in $\R_{>1}$, and the sum of the $t_i$ is $F+1$. If every $\alpha_i$ is positive,
then $\alpha^{+}=\alpha$, $\alpha^{-}=0$, $h_{a}=0$, and $E^{+}=F=E_k$ is the exponent of
Section~\ref{s:main}.

\subsection{Four lemmas}
The first lemma writes $S_{\alpha}(n)$ as a sum of denumerants with a weight.

\begin{lem}\label{lem:divisor-sum}
	Let $n\in\N_{\geq1}$. Then
	\begin{equation}\label{eq:divisor-sum}
		S_{\alpha}(n)=\sum_{a\in\N_{\geq1}^k}a_1^{\alpha_1^{-}}\dotsm a_k^{\alpha_k^{-}}\,D_{a}^{+}(n-h_{a}),
	\end{equation}
	and $D_{a}^{+}(n-h_{a})=0$ unless $d_k(a)\mid n$ and $a_i\leq n$ for every $i$; so the sum is
	finite. Moreover, for every $a\in\N_{\geq1}^k$,
	\begin{equation}\label{eq:positive-solutions}
		D_{a}^{+}(n-h_{a})=\sum_{\substack{x\in\N_{\geq1}^k\\a_1x_1+\dotsb+a_kx_k=n}}x_1^{\alpha_1^{+}}\dotsm x_k^{\alpha_k^{+}}.
	\end{equation}
\end{lem}

\begin{proof}
	\begin{itemize}
		\item[$\scriptstyle\blacklozenge$] \emph{Step~1: expanding the product.} Let $y\in\N_{\geq1}^k$ with $y_1+\dotsb+y_k=n$. By
		definition, $\sigma_{\alpha_i}(y_i)=\sum_{x_i\mid y_i}x_i^{\alpha_i}$ for every $i$. We multiply these
		$k$ sums and expand the product: we get one term for every choice of a divisor $x_1$ of $y_1$,
		a divisor $x_2$ of $y_2$, \dots, a divisor $x_k$ of $y_k$, and this term is the product of the
		chosen powers. So
		\begin{equation*}
			\sigma_{\alpha_1}(y_1)\dotsm\sigma_{\alpha_k}(y_k)
			=\sum_{x_1\mid y_1}\dotsb\sum_{x_k\mid y_k}x_1^{\alpha_1}\dotsm x_k^{\alpha_k},
		\end{equation*}
		and, summing over $y$,
		\begin{equation}\label{eq:expanded}
			S_{\alpha}(n)
			=\sum_{\substack{y\in\N_{\geq1}^k\\y_1+\dotsb+y_k=n}}\ \sum_{x_1\mid y_1}\dotsb\sum_{x_k\mid y_k}x_1^{\alpha_1}\dotsm x_k^{\alpha_k}.
		\end{equation}

		\item[$\scriptstyle\blacklozenge$] \emph{Step~2: the complementary divisors.} In a term of
		Equation~\eqref{eq:expanded}, let $a_i\coloneqq\frac{y_i}{x_i}$ for every $i$; it is a positive
		integer, since $x_i$ divides $y_i$, and $y_i=a_ix_i$. The condition $y_1+\dotsb+y_k=n$ becomes
		$a_1x_1+\dotsb+a_kx_k=n$. Conversely, given $a,x\in\N_{\geq1}^k$ with
		$a_1x_1+\dotsb+a_kx_k=n$, the numbers $y_i\coloneqq a_ix_i$ satisfy $y_1+\dotsb+y_k=n$, and $x_i$
		divides $y_i$. So the terms of Equation~\eqref{eq:expanded} correspond one to one to the pairs
		$(a,x)\in\N_{\geq1}^k\times\N_{\geq1}^k$ with $a_1x_1+\dotsb+a_kx_k=n$, and
		\begin{equation}\label{eq:pairs}
			S_{\alpha}(n)=\sum_{\substack{(a,x)\in\N_{\geq1}^k\times\N_{\geq1}^k\\a_1x_1+\dotsb+a_kx_k=n}}x_1^{\alpha_1}\dotsm x_k^{\alpha_k}.
		\end{equation}
		These pairs are finitely many, since $1\leq a_i\leq n$ and $1\leq x_i\leq n$.

		\item[$\scriptstyle\blacklozenge$] \emph{Step~3: the exchange of roles.} In Equation~\eqref{eq:pairs} the powers with a
		negative exponent $\alpha_i$ fall on $x_i$, whereas Theorem~\ref{thm:main} needs non-negative
		exponents. We move them on $a_i$. The equation $a_1x_1+\dotsb+a_kx_k=n$ contains $a_i$ and $x_i$
		only through their product $a_ix_i$, so the two numbers can exchange their roles. For a pair $(a,x)$
		of Equation~\eqref{eq:pairs} let $(a',x')$ be the pair with
		\begin{equation*}
			a_i'\coloneqq x_i,\quad x_i'\coloneqq a_i\quad\text{if }\alpha_i<0,
			\qquad\qquad
			a_i'\coloneqq a_i,\quad x_i'\coloneqq x_i\quad\text{if }\alpha_i>0.
		\end{equation*}
		Since $a_i'x_i'=a_ix_i$ for every $i$, the pair $(a',x')$ satisfies the same equation; and
		doing the exchange twice gives back $(a,x)$. So $(a,x)\mapsto(a',x')$ is a bijection of the set of
		the pairs of Equation~\eqref{eq:pairs} onto itself. Factor by factor, the term of $(a,x)$ is
		\begin{equation*}
			x_i^{\alpha_i}=x_i'^{\,\alpha_i}\quad\text{if }\alpha_i>0,
			\qquad\qquad
			x_i^{\alpha_i}=a_i'^{\,\alpha_i}\quad\text{if }\alpha_i<0 .
		\end{equation*}
		Since $(a',x')$ runs over all the pairs exactly once when $(a,x)$ does, we may sum over
		$(a',x')$; writing again $(a,x)$ in place of $(a',x')$, Equation~\eqref{eq:pairs} becomes
		\begin{equation*}
			S_{\alpha}(n)=\sum_{\substack{(a,x)\in\N_{\geq1}^k\times\N_{\geq1}^k\\a_1x_1+\dotsb+a_kx_k=n}}
			\ \prod_{i\in\{1,\dotsc,k\}\colon\alpha_i<0}a_i^{\alpha_i}
			\cdot\prod_{i\in\{1,\dotsc,k\}\colon\alpha_i>0}x_i^{\alpha_i}.
		\end{equation*}

		\item[$\scriptstyle\blacklozenge$] \emph{Step~4: grouping by $a$.} We have $\alpha_i^{-}=\alpha_i$ if $\alpha_i<0$ and
		$\alpha_i^{-}=0$ if $\alpha_i>0$; since $a_i^{0}=1$,
		\begin{equation*}
			\prod_{i\in\{1,\dotsc,k\}\colon\alpha_i<0}a_i^{\alpha_i}=a_1^{\alpha_1^{-}}\dotsm a_k^{\alpha_k^{-}},
			\qquad\text{and likewise}\qquad
			\prod_{i\in\{1,\dotsc,k\}\colon\alpha_i>0}x_i^{\alpha_i}=x_1^{\alpha_1^{+}}\dotsm x_k^{\alpha_k^{+}}.
		\end{equation*}
		The first product depends only on $a$. So, summing first over $x$ for a fixed $a$ and then over
		$a$,
		\begin{equation}\label{eq:exchange}
			S_{\alpha}(n)=\sum_{a\in\N_{\geq1}^k}a_1^{\alpha_1^{-}}\dotsm a_k^{\alpha_k^{-}}
			\sum_{\substack{x\in\N_{\geq1}^k\\a_1x_1+\dotsb+a_kx_k=n}}x_1^{\alpha_1^{+}}\dotsm x_k^{\alpha_k^{+}}.
		\end{equation}

		\item[$\scriptstyle\blacklozenge$] \emph{Step~5: the positive solutions.} Let $a\in\N_{\geq1}^k$. In the inner sum of
		Equation~\eqref{eq:exchange} we also admit $x_i=0$ for the indices with $\alpha_i>0$: such a
		solution has weight $0$, since its weight contains the factor $0^{\alpha_i^{+}}=0^{\alpha_i}=0$. So
		the inner sum is the sum of $x_1^{\alpha_1^{+}}\dotsm x_k^{\alpha_k^{+}}$ over the $x\in\N^k$
		with $a_1x_1+\dotsb+a_kx_k=n$ and $x_i\geq1$ for the indices with $\alpha_i<0$. For these
		indices let $x_i''\coloneqq x_i-1\in\N$, and for the others let $x_i''\coloneqq x_i$. We split
		the indices according to the sign of $\alpha_i$, which is never $0$; then we replace each $x_i$
		with $x_i''$ or $x_i''+1$; we expand $a_i(x_i''+1)$; and finally we join the first two sums and
		use the definition of $h_{a}$:
		\begin{align*}
			a_1x_1+\dotsb+a_kx_k
			&=\sum_{i\in\{1,\dotsc,k\}\colon\alpha_i>0}a_ix_i+\sum_{i\in\{1,\dotsc,k\}\colon\alpha_i<0}a_ix_i\\
			&=\sum_{i\in\{1,\dotsc,k\}\colon\alpha_i>0}a_ix_i''+\sum_{i\in\{1,\dotsc,k\}\colon\alpha_i<0}a_i\cdot(x_i''+1)\\
			&=\sum_{i\in\{1,\dotsc,k\}\colon\alpha_i>0}a_ix_i''+\sum_{i\in\{1,\dotsc,k\}\colon\alpha_i<0}a_ix_i''+\sum_{i\in\{1,\dotsc,k\}\colon\alpha_i<0}a_i\\
			&=a_1x_1''+\dotsb+a_kx_k''+h_{a},
		\end{align*}
		so the equation becomes $a_1x_1''+\dotsb+a_kx_k''=n-h_{a}$. The weight does not change: if
		$\alpha_i<0$, then $\alpha_i^{+}=0$ and $x_i^{\alpha_i^{+}}=1=x_i''^{\,\alpha_i^{+}}$; if $\alpha_i>0$, then
		$x_i''=x_i$. The map $x\mapsto x''$ is a bijection of these solutions onto $R_k(n-h_{a})$, with
		inverse $x_i=x_i''+1$ if $\alpha_i<0$ and $x_i=x_i''$ otherwise; if $n<h_{a}$, both sets are empty.
		Hence the inner sum is
		$\sum_{x''\in R_k(n-h_{a})}x_1''^{\,\alpha_1^{+}}\dotsm x_k''^{\,\alpha_k^{+}}=D_k^{\alpha^{+}}(n-h_{a})=D_{a}^{+}(n-h_{a})$.
		This proves Equation~\eqref{eq:positive-solutions}, and with Equation~\eqref{eq:exchange} it
		gives Equation~\eqref{eq:divisor-sum}.

		\item[$\scriptstyle\blacklozenge$] \emph{Step~6: the vanishing terms.} If $D_{a}^{+}(n-h_{a})\neq0$, by
		Equation~\eqref{eq:positive-solutions} there is $x\in\N_{\geq1}^k$ with
		$a_1x_1+\dotsb+a_kx_k=n$. Then $a_i\leq a_ix_i\leq n$ for every $i$, and $d_k(a)$, which divides
		every $a_i$, divides $n$.
	\end{itemize}
\end{proof}

So the exponents $\alpha_i<0$ do not enter the weight of the denumerant: on these coordinates the
roles of $a_i$ and $x_i$ are exchanged, and they give the factor $a_1^{\alpha_1^{-}}\dotsm a_k^{\alpha_k^{-}}$. From now on
we write the two products in the notation of Section~\ref{s:intro}: $a^{\alpha^{-}}=a_1^{\alpha_1^{-}}\dotsm a_k^{\alpha_k^{-}}$ and
$x^{\alpha^{+}}=x_1^{\alpha_1^{+}}\dotsm x_k^{\alpha_k^{+}}$. The denumerant
that remains has the weight $x^{\alpha^{+}}$, with non-negative exponents, and the next lemma bounds
it with Theorem~\ref{thm:main}.

\begin{lem}\label{lem:positive-bounds}
	Let $a\in\N_{\geq1}^k$, and let $n\in\N_{\geq1}$ be a multiple of $d_k(a)$. Then:
	\begin{enumerate}
		\item[\textup{(i)}] $D_{a}^{+}(n-h_{a})\leq C_{a}^{+}\,(n+H_{a})^{E^{+}}$;
		\item[\textup{(ii)}] if $n\geq H_{a}$, then $D_{a}^{+}(n-h_{a})\geq C_{a}^{+}\,(n-H_{a})^{E^{+}}$.
	\end{enumerate}
\end{lem}

\begin{proof}
	We use Theorem~\ref{thm:main} for the sequence $\alpha^{+}$, whose exponent is $E^{+}$ and
	whose shift is $\omega_k^{+}$. The number $n-h_{a}$ is a multiple of $d_k(a)$, since $d_k(a)$
	divides $n$ and every $a_i$.

	\emph{The shifts.} We have $\nabla_k=b_2+\frac12\sum_{i=3}^{k}b_i\leq H_{a}$. Since
		$\alpha_i^{+}=0$ exactly when $\alpha_i<0$, we have $h_{a}=\sum_{i=1}^{k}\delta_{\alpha_i^{+},0}\,a_i$.
		By the definitions of $\Delta_k$ and $\omega_k^{+}$,
		\begin{equation*}
			\Delta_k+\omega_k^{+}
			=\sum_{i=2}^{k}(b_i-a_i)+\sum_{i=2}^{k}\bigl(1-\delta_{\alpha_i^{+},0}\bigr)\,a_i
			=\sum_{i=2}^{k}b_i-\sum_{i=2}^{k}\delta_{\alpha_i^{+},0}\,a_i,
		\end{equation*}
		hence, adding $h_{a}$,
		\begin{equation}\label{eq:shift-H}
			h_{a}+\Delta_k+\omega_k^{+}=\delta_{\alpha_1^{+},0}\,a_1+\sum_{i=2}^{k}b_i\leq H_{a}.
	\end{equation}

	\emph{Part~\textup{(i)}.} If $n<h_{a}$, then $D_{a}^{+}(n-h_{a})=0$, and the right side is non-negative. If
	$n\geq h_{a}$, then $n-h_{a}\in\N$, and by Theorem~\ref{thm:main}~(i), and since $n-h_{a}\leq n$,
	$\nabla_k\leq H_{a}$ and $E^{+}\geq0$,
	\begin{equation*}
		D_{a}^{+}(n-h_{a})\leq C_{a}^{+}\,(n-h_{a}+\nabla_k)^{E^{+}}
		\leq C_{a}^{+}\,(n+H_{a})^{E^{+}}.
	\end{equation*}

	\emph{Part~\textup{(ii)}.} Let $n\geq H_{a}$. By Equation~\eqref{eq:shift-H}, $n-h_{a}\geq\Delta_k+\omega_k^{+}\geq0$.
	So Theorem~\ref{thm:main}~(ii), and then Equation~\eqref{eq:shift-H} again, give
	\begin{equation*}
		D_{a}^{+}(n-h_{a})\geq C_{a}^{+}\,\bigl(n-(h_{a}+\Delta_k+\omega_k^{+})\bigr)^{E^{+}}
		\geq C_{a}^{+}\,(n-H_{a})^{E^{+}}.\qedhere
	\end{equation*}
\end{proof}

The third lemma sums the constants of these bounds over all the vectors $a$.

\begin{lem}\label{lem:constant-sum}
	Let $n\in\N_{\geq1}$, and let
	\begin{equation*}
		c_{\alpha}\coloneqq\frac{\Gamma(\alpha_1^{+}+1)\dotsm\Gamma(\alpha_k^{+}+1)}{\Gamma(E^{+}+1)},\qquad
		\kappa_{\alpha}\coloneqq c_{\alpha}\cdot\frac{\zeta(t_1)\dotsm\zeta(t_k)}{\zeta(F+1)}.
	\end{equation*}
	Then, for every $a\in\N_{\geq1}^k$,
	\begin{align}\label{eq:C-weight}
		a_1^{\alpha_1^{-}}\dotsm a_k^{\alpha_k^{-}}\,C_{a}^{+}=c_{\alpha}\cdot\frac{d_k(a)}{a_1^{t_1}\dotsm a_k^{t_k}},
		\qquad
		\text{and}
		\qquad
		n^{E^{+}}\sum_{\substack{a\in\N_{\geq1}^k\\d_k(a)\mid n}}a_1^{\alpha_1^{-}}\dotsm a_k^{\alpha_k^{-}}\,C_{a}^{+}=\kappa_{\alpha}\,n^{E^{+}-F}\sigma_F(n).
	\end{align}
\end{lem}

\begin{proof}
	\begin{itemize}
		\item[$\scriptstyle\blacklozenge$] \emph{Step~1: the constant.} By the first formula of Lemma~\ref{lem:C-rec}, for the
		sequence $\alpha^{+}$ and the numbers $a_1,\dotsc,a_k$, and by the definition of $c_{\alpha}$,
		\begin{equation*}
			C_{a}^{+}
			=\frac{d_k(a)\,\Gamma(\alpha_1^{+}+1)\dotsm\Gamma(\alpha_k^{+}+1)}{\Gamma(E^{+}+1)\,a_1^{\alpha_1^{+}+1}\dotsm a_k^{\alpha_k^{+}+1}}
			=c_{\alpha}\cdot\frac{d_k(a)}{a_1^{\alpha_1^{+}+1}\dotsm a_k^{\alpha_k^{+}+1}}.
		\end{equation*}
		We multiply by $a_1^{\alpha_1^{-}}\dotsm a_k^{\alpha_k^{-}}$ and collect the powers of each $a_i$:
		\begin{equation*}
			a_1^{\alpha_1^{-}}\dotsm a_k^{\alpha_k^{-}}\,C_{a}^{+}
			=c_{\alpha}\,d_k(a)\cdot\frac{a_1^{\alpha_1^{-}}}{a_1^{\alpha_1^{+}+1}}\dotsm\frac{a_k^{\alpha_k^{-}}}{a_k^{\alpha_k^{+}+1}}
			=\frac{c_{\alpha}\,d_k(a)}{a_1^{\alpha_1^{+}-\alpha_1^{-}+1}\dotsm a_k^{\alpha_k^{+}-\alpha_k^{-}+1}}.
		\end{equation*}
		For every $i$ we have $\alpha_i^{+}-\alpha_i^{-}=\abs{\alpha_i}$: if $\alpha_i>0$, then
		$\alpha_i^{+}-\alpha_i^{-}=\alpha_i-0$; if $\alpha_i<0$, then $\alpha_i^{+}-\alpha_i^{-}=0-\alpha_i$. So
		the exponent of $a_i$ is $\abs{\alpha_i}+1=t_i$, and this is the first formula of Equation~\eqref{eq:C-weight}.

		\item[$\scriptstyle\blacklozenge$] \emph{Step~2: the vectors with coprime entries.} Let $v\in\N_{\geq1}^k$, let
		$g\coloneqq\gcd(v_1,\dotsc,v_k)$ and $v_i'\coloneqq\frac{v_i}{g}$ for every $i$. Then
		$v_1',\dotsc,v_k'$ are positive integers with $\gcd(v_1',\dotsc,v_k')=1$, and $v_i=g\,v_i'$.
		Conversely, given $g\in\N_{\geq1}$ and $v'\in\N_{\geq1}^k$ with $\gcd(v_1',\dotsc,v_k')=1$, the
		numbers $v_i\coloneqq g\,v_i'$ have greatest common divisor $g$. So $v\mapsto(g,v')$ is a
		bijection. All the series below have non-negative terms and converge, since every $t_i>1$; so
		we may expand the products of series and sum their terms in any order. Moreover
		\begin{equation*}
			t_1+\dotsb+t_k=\abs{\alpha_1}+\dotsb+\abs{\alpha_k}+k=\lVert\alpha\rVert+k=F+1 .
		\end{equation*}
		Hence, expanding the product, then writing $v_i=g\,v_i'$, and finally collecting the powers
		of $g$,
		\begin{align*}
			\zeta(t_1)\dotsm\zeta(t_k)
			&=\sum_{v_1\in\N_{\geq1}}\frac{1}{v_1^{t_1}}\dotsm\sum_{v_k\in\N_{\geq1}}\frac{1}{v_k^{t_k}}
			=\sum_{v\in\N_{\geq1}^k}\frac{1}{v_1^{t_1}\dotsm v_k^{t_k}}\\
			&=\sum_{g\in\N_{\geq1}}\ \sum_{\substack{v'\in\N_{\geq1}^k\\\gcd(v_1',\dotsc,v_k')=1}}\frac{1}{(g\,v_1')^{t_1}\dotsm(g\,v_k')^{t_k}}\\
			&=\sum_{g\in\N_{\geq1}}\frac{1}{g^{t_1+\dotsb+t_k}}\cdot\sum_{\substack{v'\in\N_{\geq1}^k\\\gcd(v_1',\dotsc,v_k')=1}}\frac{1}{v_1'^{\,t_1}\dotsm v_k'^{\,t_k}}
			=\zeta(F+1)\cdot\sum_{\substack{v'\in\N_{\geq1}^k\\\gcd(v_1',\dotsc,v_k')=1}}\frac{1}{v_1'^{\,t_1}\dotsm v_k'^{\,t_k}}.
		\end{align*}
		Dividing by $\zeta(F+1)$, and by the definition of $\kappa_{\alpha}$,
		\begin{equation}\label{eq:coprime-sum}
			\sum_{\substack{v'\in\N_{\geq1}^k\\\gcd(v_1',\dotsc,v_k')=1}}\frac{1}{v_1'^{\,t_1}\dotsm v_k'^{\,t_k}}
			=\frac{\zeta(t_1)\dotsm\zeta(t_k)}{\zeta(F+1)}=\frac{\kappa_{\alpha}}{c_{\alpha}}.
		\end{equation}

		\item[$\scriptstyle\blacklozenge$] \emph{Step~3: the sum.} We write every $a\in\N_{\geq1}^k$ as in Step~2:
		$a_i=g\,v_i'$, with $g=\gcd(a_1,\dotsc,a_k)=d_k(a)$ and $\gcd(v_1',\dotsc,v_k')=1$. The
		condition $d_k(a)\mid n$ becomes $g\mid n$. By the first formula of Equation~\eqref{eq:C-weight}, and since
		$t_1+\dotsb+t_k=F+1$,
		\begin{equation*}
			a_1^{\alpha_1^{-}}\dotsm a_k^{\alpha_k^{-}}\,C_{a}^{+}
			=\frac{c_{\alpha}\,g}{(g\,v_1')^{t_1}\dotsm(g\,v_k')^{t_k}}
			=\frac{c_{\alpha}\,g}{g^{F+1}\,v_1'^{\,t_1}\dotsm v_k'^{\,t_k}}
			=\frac{c_{\alpha}}{g^{F}}\cdot\frac{1}{v_1'^{\,t_1}\dotsm v_k'^{\,t_k}}.
		\end{equation*}
		We sum first over $v'$ for a fixed divisor $g$ of $n$, then over $g$; by
		Equation~\eqref{eq:coprime-sum},
		\begin{equation*}
			\sum_{\substack{a\in\N_{\geq1}^k\\d_k(a)\mid n}}a_1^{\alpha_1^{-}}\dotsm a_k^{\alpha_k^{-}}\,C_{a}^{+}
			=c_{\alpha}\sum_{g\mid n}\frac{1}{g^{F}}\cdot\sum_{\substack{v'\in\N_{\geq1}^k\\\gcd(v_1',\dotsc,v_k')=1}}\frac{1}{v_1'^{\,t_1}\dotsm v_k'^{\,t_k}}
			=\kappa_{\alpha}\sum_{g\mid n}\frac{1}{g^{F}}.
		\end{equation*}
		Finally we multiply by $n^{E^{+}}=n^{E^{+}-F}\cdot n^{F}$, and use that $g\mapsto\frac ng$ is a
		bijection of the set of the divisors of $n$:
		\begin{equation*}
			n^{E^{+}}\sum_{g\mid n}\frac{1}{g^{F}}
			=n^{E^{+}-F}\sum_{g\mid n}\Bigl(\frac{n}{g}\Bigr)^{F}
			=n^{E^{+}-F}\sum_{g\mid n}g^{F}
			=n^{E^{+}-F}\sigma_F(n).\qedhere
		\end{equation*}
	\end{itemize}
\end{proof}

The last lemma collects three rough estimates, used for the vectors $a$ with large entries.

\begin{lem}\label{lem:sizes}
	Let $a\in\N_{\geq1}^k$. Then:
	\begin{enumerate}
		\item[\textup{(i)}] $d_k(a)H_{a}\leq k\,a_1\dotsm a_k$ and $H_{a}\leq k\,a_1\dotsm a_k$;
		\item[\textup{(ii)}] $d_k(a)\leq a_1^{\abs{\alpha_1}/\lVert\alpha\rVert}\dotsm a_k^{\abs{\alpha_k}/\lVert\alpha\rVert}$;
		\item[\textup{(iii)}] for every $n\in\N_{\geq1}$,
		$\displaystyle a^{\alpha^{-}}D_{a}^{+}(n-h_{a})\leq n^{E^{+}}\prod_{i=1}^{k}a_i^{\abs{\alpha_i}/\lVert\alpha\rVert-t_i}$.
	\end{enumerate}
\end{lem}

\begin{proof}
	\emph{Part~\textup{(i)}.} Since $d_k(a)$ divides $a_2$, we have $d_k(a)a_1\leq a_1a_2\leq a_1\dotsm a_k$, the
	last inequality because the other factors are at least $1$. Let $i\in\{2,\dotsc,k\}$. As in the
	proof of Lemma~\ref{lem:C-rec}, $b_i=\frac{d_{i-1}a_i}{d_i}$. Moreover $d_{i-1}\leq a_1$, since
	$d_{i-1}$ divides $a_1$, and $d_i\geq d_k(a)$, since $d_k$ divides $d_i$. Hence
	$d_k(a)b_i\leq a_1a_i\leq a_1\dotsm a_k$. Summing these $k$ inequalities we obtain the first
	inequality, and the second follows since $d_k(a)\geq1$.

	\emph{Part~\textup{(ii)}.} The exponents $\frac{\abs{\alpha_i}}{\lVert\alpha\rVert}$ are non-negative and their
	sum is $1$. The number $d_k(a)$ is at most $a_i$ for every $i$; so $d_k(a)$, which is the product of
	its powers with these exponents, is at most the product of the corresponding powers of the $a_i$.

	\emph{Part~\textup{(iii)}.} By Equation~\eqref{eq:positive-solutions}, $D_{a}^{+}(n-h_{a})$ is the sum of
	$x^{\alpha^{+}}$ over the $x\in\N_{\geq1}^k$ with $a_1x_1+\dotsb+a_kx_k=n$. Let
	$j\in\{1,\dotsc,k\}$. For such an $x$ we have $a_ix_i\leq n$, that is,
	$1\leq x_i\leq\frac{n}{a_i}$, for every $i$; and $x$ is determined by the numbers $x_i$ with
	$i\neq j$, since $x_j=\frac{1}{a_j}\bigl(n-\sum_{i\neq j}a_ix_i\bigr)$. So there are at most
	$\prod_{i\neq j}\frac{n}{a_i}$ such $x$, and each weight $x^{\alpha^{+}}$ is at most
	$\prod_{i=1}^{k}\bigl(\frac{n}{a_i}\bigr)^{\alpha_i^{+}}$. Hence, since
	$\abs{\alpha^{+}}+k-1=E^{+}$,
	\begin{equation*}
		D_{a}^{+}(n-h_{a})
		\leq\prod_{i\neq j}\frac{n}{a_i}\cdot\prod_{i=1}^{k}\Bigl(\frac{n}{a_i}\Bigr)^{\alpha_i^{+}}
		=\frac{n^{E^{+}}\,a_j}{a_1^{\alpha_1^{+}+1}\dotsm a_k^{\alpha_k^{+}+1}}.
	\end{equation*}
	This holds for every $j$. We raise the bound for the index $j$ to the power
	$\frac{\abs{\alpha_j}}{\lVert\alpha\rVert}$ and multiply over $j$; since these exponents are
	non-negative with sum $1$, the left side becomes $D_{a}^{+}(n-h_{a})$ again, and
	\begin{equation*}
		D_{a}^{+}(n-h_{a})\leq\frac{n^{E^{+}}\,a_1^{\abs{\alpha_1}/\lVert\alpha\rVert}\dotsm a_k^{\abs{\alpha_k}/\lVert\alpha\rVert}}{a_1^{\alpha_1^{+}+1}\dotsm a_k^{\alpha_k^{+}+1}}.
	\end{equation*}
	Finally, as in Step~1 of the proof of Lemma~\ref{lem:constant-sum}, multiplying by
	$a^{\alpha^{-}}$ turns the denominator into $a_1^{t_1}\dotsm a_k^{t_k}$.
\end{proof}

\subsection{The main term}
\begin{thm}\label{thm:divisor}
	Let $k\in\N_{\geq2}$ and $\alpha\in(\R\setminus\{0\})^k$ with $\lVert\alpha\rVert\in\R_{>1}$. Recall
	from Section~\ref{ss:setting} and Lemma~\ref{lem:constant-sum} that
	\begin{gather*}
		E^{+}=\alpha_1^{+}+\dotsb+\alpha_k^{+}+k-1,\qquad
		F=\abs{\alpha_1}+\dotsb+\abs{\alpha_k}+k-1,\qquad
		t_i=\abs{\alpha_i}+1,\\
		c_{\alpha}=\frac{\Gamma(\alpha_1^{+}+1)\dotsm\Gamma(\alpha_k^{+}+1)}{\Gamma(E^{+}+1)},\qquad
		\kappa_{\alpha}=c_{\alpha}\cdot\frac{\zeta(t_1)\dotsm\zeta(t_k)}{\zeta(F+1)}.
	\end{gather*}
	For $\theta\in(0,1)$
	and $i\in\{1,\dotsc,k\}$ let
	\begin{equation*}
		s_i(\theta)\coloneqq1+\abs{\alpha_i}-\theta-(1-\theta)\,\frac{\abs{\alpha_i}}{\lVert\alpha\rVert},\qquad
		s_i'(\theta)\coloneqq1+\abs{\alpha_i}-\theta-\frac{\abs{\alpha_i}}{\lVert\alpha\rVert}.
	\end{equation*}
	Let $n\in\N_{\geq1}$.
	\begin{enumerate}
		\item[\textup{(i)}] For every $\theta\in(0,1)$ with $s_i(\theta)\in\R_{>1}$ for every $i$,
		\begin{equation*}
			S_{\alpha}(n)\geq\kappa_{\alpha}\,n^{E^{+}-F}\sigma_F(n)
			-(E^{+}+1)\,c_{\alpha}\,k^{\theta}\,\zeta\bigl(s_1(\theta)\bigr)\dotsm\zeta\bigl(s_k(\theta)\bigr)\,n^{E^{+}-\theta}.
		\end{equation*}
		\item[\textup{(ii)}] For every $\theta\in(0,1)$ with $s_i'(\theta)\in\R_{>1}$ for every $i$,
		\begin{equation*}
			S_{\alpha}(n)\leq\kappa_{\alpha}\,n^{E^{+}-F}\sigma_F(n)
			+k^{\theta}\Bigl(E^{+}\,2^{E^{+}-1}c_{\alpha}\prod_{i=1}^{k}\zeta\bigl(s_i(\theta)\bigr)
			+\prod_{i=1}^{k}\zeta\bigl(s_i'(\theta)\bigr)\Bigr)\,n^{E^{+}-\theta}.
		\end{equation*}
	\end{enumerate}
	In particular, $\lim_{n\to\infty}\frac{S_{\alpha}(n)}{n^{E^{+}-F}\sigma_F(n)}=\kappa_{\alpha}$.
\end{thm}

Both conditions on $\theta$ hold for every small $\theta$. Indeed $s_i'(\theta)>1$ means
$\theta<\abs{\alpha_i}\bigl(1-\frac{1}{\lVert\alpha\rVert}\bigr)$, and the right side is positive
since $\lVert\alpha\rVert>1$; moreover $s_i(\theta)=s_i'(\theta)+\theta\,\frac{\abs{\alpha_i}}{\lVert\alpha\rVert}\geq s_i'(\theta)$.
If every $\alpha_i$ is positive, then $E^{+}=F$ and the main term is $\kappa_{\alpha}\,\sigma_{E^{+}}(n)$.

\begin{proof}
	We have $E^{+}\geq k-1\geq1$. All the sums over $a$ in this proof are over $a\in\N_{\geq1}^k$.

	\begin{itemize}
		\item[$\scriptstyle\blacklozenge$] \emph{Step~1: a convergent series.} Let $\theta\in(0,1)$ with $s_i(\theta)>1$ for every $i$. We
		prove that
		\begin{equation}\label{eq:Z-theta}
			T_{\theta}\coloneqq\sum_{a}a^{\alpha^{-}}\,C_{a}^{+}H_{a}^{\theta}
			\leq c_{\alpha}\,k^{\theta}\,\zeta\bigl(s_1(\theta)\bigr)\dotsm\zeta\bigl(s_k(\theta)\bigr).
		\end{equation}
		By the first formula of Equation~\eqref{eq:C-weight}, then by Lemma~\ref{lem:sizes}~(i) and~(ii),
		\begin{align*}
			a^{\alpha^{-}}\,C_{a}^{+}H_{a}^{\theta}
			&=c_{\alpha}\,\frac{d_k(a)^{1-\theta}\,(d_k(a)H_{a})^{\theta}}{a_1^{t_1}\dotsm a_k^{t_k}}\\
			&\leq c_{\alpha}\,k^{\theta}\prod_{i=1}^{k}a_i^{(1-\theta)\abs{\alpha_i}/\lVert\alpha\rVert+\theta-t_i}
			=c_{\alpha}\,k^{\theta}\prod_{i=1}^{k}a_i^{-s_i(\theta)}.
		\end{align*}
		Summing over $a$, the sum of the product is
		$\zeta\bigl(s_1(\theta)\bigr)\dotsm\zeta\bigl(s_k(\theta)\bigr)$, which is finite since every
		$s_i(\theta)>1$. This is Equation~\eqref{eq:Z-theta}.

		\item[$\scriptstyle\blacklozenge$] \emph{Step~2: two elementary inequalities.} Let $\theta$ be as in Step~1, so that $T_{\theta}$ is
		finite. If $H_{a}>n$, then
		$1<\bigl(\frac{H_{a}}{n}\bigr)^{\theta}$; if $H_{a}\leq n$, then
		$H_{a}=H_{a}^{\theta}H_{a}^{1-\theta}\leq H_{a}^{\theta}n^{1-\theta}$. Hence
		\begin{equation}\label{eq:two-tails}
			n^{E^{+}}\sum_{\substack{d_k(a)\mid n\\H_{a}>n}}a^{\alpha^{-}}\,C_{a}^{+}\leq n^{E^{+}-\theta}T_{\theta},
			\qquad
			n^{E^{+}-1}\sum_{\substack{d_k(a)\mid n\\H_{a}\leq n}}a^{\alpha^{-}}\,C_{a}^{+}H_{a}\leq n^{E^{+}-\theta}T_{\theta}.
		\end{equation}

		\item[$\scriptstyle\blacklozenge$] \emph{Step~3: the lower bound.} Let $\theta$ be as in~(i), and let $\mathcal{A}_{\leq}$ be the set
		of the $a$ with $d_k(a)\mid n$ and $H_{a}\leq n$. By Lemma~\ref{lem:divisor-sum}, and since its
		terms are non-negative,
		\begin{equation*}
			S_{\alpha}(n)\geq\sum_{a\in\mathcal{A}_{\leq}}a^{\alpha^{-}}\,D_{a}^{+}(n-h_{a}).
		\end{equation*}
		For $a\in\mathcal{A}_{\leq}$, Lemma~\ref{lem:positive-bounds}~(ii) gives
		$D_{a}^{+}(n-h_{a})\geq C_{a}^{+}(n-H_{a})^{E^{+}}$. Since $E^{+}\geq1$, we have $(1-u)^{E^{+}}\geq1-E^{+}u$ for
		$u\in[0,1]$, so $(n-H_{a})^{E^{+}}\geq n^{E^{+}}-E^{+}\,H_{a}\,n^{E^{+}-1}$. Hence, by Lemma~\ref{lem:constant-sum}
		and by Equation~\eqref{eq:two-tails},
		\begin{align*}
			S_{\alpha}(n)
			&\geq n^{E^{+}}\sum_{a\in\mathcal{A}_{\leq}}a^{\alpha^{-}}\,C_{a}^{+}
			-E^{+}\,n^{E^{+}-1}\sum_{a\in\mathcal{A}_{\leq}}a^{\alpha^{-}}\,C_{a}^{+}H_{a}\\
			&=\kappa_{\alpha}\,n^{E^{+}-F}\sigma_F(n)-n^{E^{+}}\sum_{\substack{d_k(a)\mid n\\H_{a}>n}}a^{\alpha^{-}}\,C_{a}^{+}
			-E^{+}\,n^{E^{+}-1}\sum_{a\in\mathcal{A}_{\leq}}a^{\alpha^{-}}\,C_{a}^{+}H_{a}\\
			&\geq\kappa_{\alpha}\,n^{E^{+}-F}\sigma_F(n)-(E^{+}+1)\,n^{E^{+}-\theta}T_{\theta},
		\end{align*}
		and Equation~\eqref{eq:Z-theta} gives~(i).

		\item[$\scriptstyle\blacklozenge$] \emph{Step~4: the upper bound, small $H_{a}$.} Let $\theta$ be as in~(ii); then also
		$s_i(\theta)\geq s_i'(\theta)>1$, so Equation~\eqref{eq:Z-theta} holds. Let $\mathcal{A}_{\leq}$
		be as in Step~3. For $a\in\mathcal{A}_{\leq}$, Lemma~\ref{lem:positive-bounds}~(i) gives
		$D_{a}^{+}(n-h_{a})\leq C_{a}^{+}(n+H_{a})^{E^{+}}$. By the mean value theorem,
		$(n+H_{a})^{E^{+}}-n^{E^{+}}=E^{+}\,(n+\xi)^{E^{+}-1}H_{a}$ for some $\xi\in(0,H_{a})$; since $E^{+}-1\geq0$ and
		$n+\xi\leq2n$, we get $(n+H_{a})^{E^{+}}\leq n^{E^{+}}+E^{+}\,2^{E^{+}-1}n^{E^{+}-1}H_{a}$. Hence, since the terms are
		non-negative, by Lemma~\ref{lem:constant-sum}, by Equation~\eqref{eq:two-tails} and by
		Equation~\eqref{eq:Z-theta},
		\begin{align*}
			\sum_{a\in\mathcal{A}_{\leq}}a^{\alpha^{-}}\,D_{a}^{+}(n-h_{a})
			&\leq\kappa_{\alpha}\,n^{E^{+}-F}\sigma_F(n)+E^{+}\,2^{E^{+}-1}n^{E^{+}-\theta}T_{\theta}\\
			&\leq\kappa_{\alpha}\,n^{E^{+}-F}\sigma_F(n)+E^{+}\,2^{E^{+}-1}c_{\alpha}\,k^{\theta}\,\zeta\bigl(s_1(\theta)\bigr)\dotsm\zeta\bigl(s_k(\theta)\bigr)\,n^{E^{+}-\theta}.
		\end{align*}

		\item[$\scriptstyle\blacklozenge$] \emph{Step~5: the upper bound, large $H_{a}$.} Let $\mathcal{A}_{>}$ be the set of the $a$ with
		$d_k(a)\mid n$ and $H_{a}>n$. For $a\in\mathcal{A}_{>}$, by Lemma~\ref{lem:sizes}~(iii), then since
		$1<\bigl(\frac{H_{a}}{n}\bigr)^{\theta}$, and finally by Lemma~\ref{lem:sizes}~(i),
		\begin{align*}
			a^{\alpha^{-}}\,D_{a}^{+}(n-h_{a})
			&\leq n^{E^{+}}\prod_{i=1}^{k}a_i^{\abs{\alpha_i}/\lVert\alpha\rVert-t_i}
			\leq n^{E^{+}-\theta}H_{a}^{\theta}\prod_{i=1}^{k}a_i^{\abs{\alpha_i}/\lVert\alpha\rVert-t_i}\\
			&\leq k^{\theta}\,n^{E^{+}-\theta}\prod_{i=1}^{k}a_i^{\abs{\alpha_i}/\lVert\alpha\rVert+\theta-t_i}
			=k^{\theta}\,n^{E^{+}-\theta}\prod_{i=1}^{k}a_i^{-s_i'(\theta)}.
		\end{align*}
		Since every $s_i'(\theta)>1$, summing over $a$ we obtain
		\begin{equation*}
			\sum_{a\in\mathcal{A}_{>}}a^{\alpha^{-}}\,D_{a}^{+}(n-h_{a})\leq k^{\theta}\,\zeta\bigl(s_1'(\theta)\bigr)\dotsm\zeta\bigl(s_k'(\theta)\bigr)\,n^{E^{+}-\theta}.
		\end{equation*}

		\item[$\scriptstyle\blacklozenge$] \emph{Step~6: conclusion.} By Lemma~\ref{lem:divisor-sum}, only the $a$ with
		$d_k(a)\mid n$ contribute, and these are in $\mathcal{A}_{\leq}$ or in $\mathcal{A}_{>}$; so Steps~4
		and~5 give~(ii). Finally, we fix a number $\theta$ that satisfies the condition of~(ii), and hence also that of~(i);
		it exists by the remark after the statement. By~(i) and~(ii), the difference between $S_{\alpha}(n)$ and
		$\kappa_{\alpha}n^{E^{+}-F}\sigma_F(n)$ is at most a constant, which does not depend on $n$, times $n^{E^{+}-\theta}$, while
		\begin{equation*}
			n^{E^{+}-F}\sigma_F(n)=n^{E^{+}}\sum_{g\mid n}\Bigl(\frac{g}{n}\Bigr)^{F}\geq n^{E^{+}},
		\end{equation*}
		since the term $g=n$ is $1$; so the quotient $\frac{S_{\alpha}(n)}{n^{E^{+}-F}\sigma_F(n)}$ differs from $\kappa_{\alpha}$ by at
		most a constant times $n^{-\theta}$, which tends to $0$.
	\end{itemize}
\end{proof}

\begin{rem}\label{rem:divisor}
	\begin{enumerate}
		\item[(a)] \emph{Equal positive exponents.} Let $u\in\R_{>0}$ with $ku>1$, and let
		$\alpha=(u,\dotsc,u)$. Then $\alpha^{+}=\alpha$, $\lVert\alpha\rVert=ku$, $t_i=u+1$ for every $i$, and
		$E^{+}=F=ku+k-1$; so
		\begin{equation*}
			c_{\alpha}=\frac{\Gamma(u+1)^{k}}{\Gamma(ku+k)},\qquad
			\kappa_{\alpha}=\frac{\Gamma(u+1)^{k}}{\Gamma(ku+k)}\cdot\frac{\zeta(u+1)^{k}}{\zeta(ku+k)}.
		\end{equation*}
		Since $E^{+}-F=0$, Theorem~\ref{thm:divisor} gives
		\begin{align*}
			\sum_{\substack{y\in\N_{\geq1}^k\\y_1+\dotsb+y_k=n}}\sigma_{u}(y_1)\dotsm\sigma_{u}(y_k)
			&\sim\frac{\Gamma(u+1)^{k}}{\Gamma(ku+k)}\cdot\frac{\zeta(u+1)^{k}}{\zeta(ku+k)}\,\sigma_{ku+k-1}(n)\\
			&=\frac{\Gamma(u+1)^{k}}{\Gamma(ku+k)}\cdot\frac{\zeta(u+1)^{k}}{\zeta(ku+k)}\,n^{ku+k-1}\,\sigma_{-(ku+k-1)}(n).
		\end{align*}
		The last equality holds since $\sigma_{v}(n)=n^{v}\sigma_{-v}(n)$ for every $v\in\R$, because
		$d\mapsto\frac nd$ is a bijection of the set of the divisors of $n$:
		$\sigma_{v}(n)=\sum_{d\mid n}\bigl(\frac nd\bigr)^{v}=n^{v}\sum_{d\mid n}d^{-v}=n^{v}\sigma_{-v}(n)$; here
		$v=ku+k-1$. These formulas hold for every real $u>0$ with $ku>1$, and for this reason they
		are written with the Gamma function. If we restrict to $u\in\N_{\geq1}$ (then $ku>1$ holds,
		since $k\geq2$), then $u+1$ and $ku+k$ are positive integers, so $\Gamma(u+1)=u!$ and
		$\Gamma(ku+k)=(ku+k-1)!$, and the formula reads
		\begin{align*}
			\sum_{\substack{y\in\N_{\geq1}^k\\y_1+\dotsb+y_k=n}}\sigma_{u}(y_1)\dotsm\sigma_{u}(y_k)
			&\sim\frac{(u!)^{k}}{(ku+k-1)!}\cdot\frac{\zeta(u+1)^{k}}{\zeta(ku+k)}\,\sigma_{ku+k-1}(n)\\
			&=\frac{(u!)^{k}}{(ku+k-1)!}\cdot\frac{\zeta(u+1)^{k}}{\zeta(ku+k)}\,n^{ku+k-1}\,\sigma_{-(ku+k-1)}(n).
		\end{align*}
		For $u=1$, since $1!=1$, the constant is
		$\kappa_{\alpha}=\frac{\zeta(2)^{k}}{(2k-1)!\,\zeta(2k)}$. For $k=2$ and $u=1$ this is
		$\kappa_{\alpha}=\frac{\zeta(2)^2}{6\,\zeta(4)}=\frac{5}{12}$, the coefficient of $\sigma_3(n)$ in
		the classical identity
		\begin{equation*}
			\sum_{y=1}^{n-1}\sigma_1(y)\,\sigma_1(n-y)=\frac{5\,\sigma_3(n)-(6n-1)\,\sigma_1(n)}{12};
		\end{equation*}
		for $k=2$ the constant $\kappa_{\alpha}$ is the one of Ramanujan~\cite{ramanujan}.
		\item[(b)] \emph{Equal negative exponents.} Let $u\in\R_{>0}$ with $ku>1$, and let
		$\alpha=(-u,\dotsc,-u)$. Then $\alpha^{+}=(0,\dotsc,0)$, $\lVert\alpha\rVert=ku$, $t_i=u+1$ for
		every $i$, $E^{+}=k-1$ and $F=ku+k-1$; so $\Gamma(E^{+}+1)=\Gamma(k)=(k-1)!$ and
		\begin{equation*}
			c_{\alpha}=\frac{1}{(k-1)!},\qquad
			\kappa_{\alpha}=\frac{1}{(k-1)!}\cdot\frac{\zeta(u+1)^{k}}{\zeta(ku+k)}.
		\end{equation*}
		Since $E^{+}-F=-ku$, Theorem~\ref{thm:divisor} gives
		\begin{align*}
			\sum_{\substack{y\in\N_{\geq1}^k\\y_1+\dotsb+y_k=n}}\sigma_{-u}(y_1)\dotsm\sigma_{-u}(y_k)
			&\sim\frac{1}{(k-1)!}\cdot\frac{\zeta(u+1)^{k}}{\zeta(ku+k)}\,n^{-ku}\,\sigma_{ku+k-1}(n)\\
			&=\frac{n^{k-1}}{(k-1)!}\cdot\frac{\zeta(u+1)^{k}}{\zeta(ku+k)}\,\sigma_{-(ku+k-1)}(n).
		\end{align*}
		The last equality holds as in~(a), since $\sigma_{ku+k-1}(n)=n^{ku+k-1}\sigma_{-(ku+k-1)}(n)$ and
		$n^{-ku}\cdot n^{ku+k-1}=n^{k-1}$.
		For $k=2$ and $u=1$ this is $\kappa_{\alpha}=\frac{\zeta(2)^2}{\zeta(4)}=\frac52$, and the sum is
		asymptotic to $\frac52\,n^{-2}\sigma_3(n)$. The constants of~(a) and~(b) have the same factor
		$\frac{\zeta(u+1)^{k}}{\zeta(ku+k)}$, which comes from the sum over the vectors $a$
		(Lemma~\ref{lem:constant-sum}); they differ only in $c_{\alpha}$. When every $\alpha_i$ is
		negative, as here, $\alpha^{+}=0$ and the proof uses Theorem~\ref{thm:main} only for $\alpha=0$,
		that is, the bounds of~\cite{ballesio-I} for the denumerant.
		\item[(c)] \emph{The hypotheses.} An exponent $\alpha_i=0$ is excluded, since then
		$t_i=1$ and the series $\zeta(t_i)$ of Lemma~\ref{lem:constant-sum} diverges. The hypothesis
		$\lVert\alpha\rVert>1$ is needed for the existence of $\theta$: as $\theta\to0$, both conditions
		$s_i(\theta)>1$ and $s_i'(\theta)>1$ become
		$\abs{\alpha_i}\bigl(1-\frac{1}{\lVert\alpha\rVert}\bigr)>0$. It comes from the rough bounds of
		Lemma~\ref{lem:sizes}~(ii) and~(iii); for $k=2$ the cases with $\lVert\alpha\rVert\leq1$ and positive
		exponents are covered by Ingham's theorem~\cite{ingham}.
		\item[(d)] \emph{The size of the error.} The constants of Theorem~\ref{thm:divisor} are not
		optimized. Numerically, for $\alpha=(1,1)$ the lower bound of Step~3, before the use of
		$\theta$, differs from $\kappa_{\alpha}\sigma_3(n)$ by at most
		$0.36\,n^{2}(1+\log n)^2$ for $2\leq n\leq400$, which suggests an error of order
		$n^{E^{+}-1}\log^2n$ for this method.
	\end{enumerate}
\end{rem}

\bigskip
\noindent\textbf{Acknowledgements.} The author used an AI assistant (Claude, Anthropic) to revise
the English wording and the organization of the text, and to check the formulas and the
numerical examples. All results and proofs are due to the author, who takes full responsibility
for the content of the paper.

\bibliographystyle{aomalpha}
\bibliography{denumerant_bounds_III}

\end{document}